\documentclass[11pt]{amsart}

\usepackage{amsxtra,amssymb,amsmath,amscd,url,listings,stmaryrd}
\usepackage[alphabetic,initials]{amsrefs}
\usepackage{scalerel,stackengine}
\stackMath
\usepackage[utf8]{inputenc}
\usepackage{eucal}
\usepackage{fullpage}
\usepackage{scrtime}
\usepackage[colorlinks]{hyperref}
\usepackage{xcolor}
\usepackage{datetime}
\usepackage{graphicx}

\shortdate

\renewcommand{\leq}{\leqslant}
\renewcommand{\geq}{\geqslant}
\newcommand\widecheck[1]{%
\savestack{\tmpbox}{\stretchto{%
  \scaleto{%
    \scalerel*[\widthof{\ensuremath{#1}}]{\kern-.4pt\bigwedge\kern-.4pt}%
    {\rule[-\textheight/2]{1ex}{\textheight}}
  }{\textheight}%
}{0.5ex}}%
\stackon[2pt]{#1}{\scalebox{-1}{\tmpbox}}%
}

\numberwithin{equation}{section}

\def\stacksum#1#2{{\stackrel{{\scriptstyle #1}}
{{\scriptstyle #2}}}}

\newcommand{\FT}{\mathrm{FT}}

\newcommand{\sym}{\mathrm{sym}}

\newcommand{\Cc}{\mathbf{C}}

\newcommand{\Aa}{\mathbf{A}}

\newcommand{\Zz}{\mathbf{Z}}

\newcommand{\Rr}{\mathbf{R}}

\newcommand{\Qq}{\mathbf{Q}}

\newcommand{\Fq}{{\mathbf{F}_q}}

\newcommand{\Fqt}{{\mathbf{F}^\times_q}}
\newcommand{\Fqqt}{{\mathbf{F}^\times_{q_0}}}

\newcommand{\Ff}{\mathbf{F}}

\newcommand{\mcW}{\mathcal{W}}

\newcommand{\KL}{\mathcal{K}\ell}

\newcommand{\mods}[1]{\,(\mathrm{mod}\,{#1})}

\newcommand{\what}{\widehat}

\newcommand{\ra}{\rightarrow}

\DeclareMathOperator{\Kl}{\mathrm{Kl}}

\newcommand{\eps}{\varepsilon}
\renewcommand{\rho}{\varrho}

\DeclareMathOperator{\GL}{GL}

\DeclareMathSymbol{\gena}{\mathord}{letters}{"3C}
\DeclareMathSymbol{\genb}{\mathord}{letters}{"3E}

\newcommand{\Xcv}{X_{\mathrm{cv}}}
\newcommand{\Xap}{X_{\mathrm{ap}}}

\theoremstyle{plain}
\newtheorem{theorem}{Theorem}[section]
\newtheorem*{theorem*}{Theorem}
\newtheorem{lemma}[theorem]{Lemma}

\newtheorem{corollary}[theorem]{Corollary}

\newtheorem{proposition}[theorem]{Proposition}

\newtheorem*{notation}{Notation}

\theoremstyle{remark}

\theoremstyle{definition}

\newtheorem{remark}[theorem]{Remark}

\newcommand{\mcC}{\mathcal{C}}

\newcommand{\rmP}{\mathrm{P}}

\newcommand{\mcF}{\mathcal{F}}

\newcommand{\mcK}{\mathcal{K}}

\newcommand{\mcJ}{\mathcal{J}}

\newcommand{\lf}{\lambda_f}

\newcommand{\vphi}{\varphi}

\renewcommand{\geq}{\geqslant}
\renewcommand{\leq}{\leqslant}
\renewcommand{\Re}{\mathfrak{Re}\,}

\newcommand{\ov}[1]{\overline{#1}}

\newcommand\sumsum{\mathop{\sum\sum}\limits}

\newcommand{\sumstar}{\sideset{}{^\star}\sum}
\newcommand{\sumsumstar}{\sideset{}{^\star}\sumsum}

\begin{document}
\title{On algebraic twists with composite moduli}

\author{Yongxiao Lin}
\address{Data Science Institute, Shandong University, Jinan 250100, China}
\email{yongxiao.lin@sdu.edu.cn}

\author{Philippe Michel}
\address{EPFL/MATH/TAN, Station 8, CH-1015 Lausanne, Switzerland }
\email{philippe.michel@epfl.ch}


\begin{abstract} 
  We study bounds for algebraic twists sums of automorphic coefficients by trace functions of composite moduli.
\end{abstract}

\thanks{Ph.\ M. was partially supported by the SNF (grant 200021\_197045).  \today\ \currenttime}

\maketitle

\hfill{\em In memory of Chandra Sekhar Raju}
\section{Introduction}
In a series of papers \cite{FKM1,KLMS,LMS} we studied the absence of correlations between the coefficients of certain automorphic $L$-functions and trace functions of prime moduli. 

More precisely, given $q$ a prime number, let
$$K:\Fq=\Zz/q\Zz\mapsto \Cc$$ be the trace function associated to a suitable $\ell$-adic middle extension sheaf $\mcF$  on the affine line $\Aa^1_{\Fq}$, geometrically irreducible and pure of weight $0$; this implies in particular that the supnorm of $K$ satisfies
$$\|K\|_\infty\leq C(\mcF)$$
where $C(\mcF)$ denotes the analytic conductor of $\mcF$, a numerical invariant attached to the Galois representation underlying $\mcF$. We now view $K$ as a $q$-periodic function on $\Zz$ via the obvious projection.

Let
$$L(\pi,s)=\sum_{n\geq 1}\frac{\lambda_\pi(n)}{n^s}=\prod_p L(\pi_p,s),\ \Re s>1$$
be an automorphic $L$-function of some degree $d\geq 2$ (normalized so that $\Re s=1/2$ is the critical line). For $V$ a smooth, compactly supported function on $\Rr_{>0}$, we consider the problem of obtaining non-trivial bounds for the correlation sums
\begin{equation}\label{nontrivialgeneric}
S_V(K;X)=\sum_{n\geq 1}\lambda_\pi(n)K(n)V\left(\frac{n}X\right)\ll X^{1-\eta}\hbox{  as $q,X\ra\infty$;}	
\end{equation}
here $\eta>0$ is some positive constant and the above bound depends implicitly on $\pi$, $C(\mcF)$ and $V$ and for $X$ varying over a range as small as possible compared to $q$. 

\subsubsection*{The convexity range} Under relatively mild conditions on $\mcF$ it is not too difficult to obtain non-trivial bounds like \eqref{nontrivialgeneric} as long as
$$X\geq q^{d/2+\delta}$$
for some $\delta>0$ (with $\eta$ depending on $\delta$) and so the first challenge is to pass the so-called the {\em convexity range}
\begin{equation}\label{convexrange}
	\Xcv:= q^{d/2}.
	\end{equation}
 Indeed passing this range for $K=\chi\bmod  q$ a non-trivial Dirichlet character, enables one to solve the subconvexity problem  for the twisted $L$-function $L(\pi\times\chi,s)$ for $\Re s=1/2$ in the large $q$-aspect.

In the three papers mentioned above, a non-trivial bound \eqref{nontrivialgeneric} was obtained for $X$ at and below the convexity range. Specifically
\begin{itemize}
\item \cite{FKM1} considered the situation where $L(\pi,s)$ is the standard $L$-function of a $\GL_{2,\Qq}$ automorphic representation (the $L$-function of a Hecke eigenform) and obtained (under some mild assumptions on $\mcF$ that are recalled below) \eqref{nontrivialgeneric} as long as
\begin{equation}\label{shortestd2}X\geq q^{1-1/4+\delta},\ \delta>0.	
\end{equation}

\item \cite{KLMS} considered the situation where $L(\pi,s)$ is the standard $L$-function of a $\GL_{3,\Qq}$ automorphic representation (of level $1$) and obtained (again under some suitable assumptions on $\mcF$)  \eqref{nontrivialgeneric} as long as
\begin{equation}\label{shortestd3}
	X\geq q^{3/2-1/6+\delta},\ \delta>0.
\end{equation}

\item \cite{LMS} considered the situation where $L(\pi,s)$ is the  Rankin--Selberg $L$-function attached to a pair $(\vphi,f)$ of
$\GL_{3,\Qq}$ and $\GL_{2,\Qq}$ automorphic forms (both of level $1$). More precisely  $\lambda_\pi$ is given by
\begin{equation}
	\label{RS23def}\lambda_\pi(n)=\sum_{mr^2=n}\lambda_\vphi(m,r)\lambda_f(m).
\end{equation} In that case \eqref{nontrivialgeneric} can be obtained for $K$ a trace function  associated with a suitably ``good" sheaf $\mcK$ (see \cite[\S 1]{LMS} for the definition of the goodness) as long as
\begin{equation}\label{shortestd6}X\geq q^{3-1/4+\delta},\ \delta>0.	
\end{equation}
\end{itemize}

\subsubsection*{The arithmetic progression range}
It is of course  desirable and often interesting to try to obtain \eqref{nontrivialgeneric} for even shorter ranges. Such a range arises when studying the distribution of  $(\lambda_\pi(n))_{n\leq X}$ in large arithmetic progressions. Given $q$ a modulus and a primitive congruence class $a\bmod  q$ (i.e., $(a,q)=1$), the goal is to improve the trivial estimate\footnote{This may require the Ramanujan--Petersson conjecture.} for the sum
\begin{equation}\label{deltasum}
\sum_{n\equiv a\bmod  q}\lambda_\pi(n)V\left(\frac{n}X\right)\ll X^{o(1)}\frac{X}q.
\end{equation}
Expressing the congruence $n\equiv a\bmod  q$ in terms of Dirichlet characters $\bmod  q$ and using the functional equation for $L(\pi\times\chi,s)$ transforms the left hand side of \eqref{deltasum} into a sum  of the shape (possibly up to some main terms)
\begin{equation}\label{deltatransform}
	\frac{X}{q^{\frac{d+1}2}}\sum_{n}\lambda_{\ov\pi}(n)\Kl_d(an;q)\widecheck{V}\left(\frac{nX}{q^d}\right)
\end{equation}
 where $\widecheck{V}(x)$ is a suitable integral transform of $V$ (depending on $d$ and the Gamma factors of $\pi$) and is rapidly decreasing and
$$\Kl_d(n; q)=\frac{1}{q^{\frac{d-1}{2}}}\sumsum_\stacksum{x_1,\cdots,x_d\in(\Zz/q\Zz)^\times}{x_1.\cdots.x_d=n}e\left(\frac{x_1+\cdots+x_d}{q}\right)$$
denotes the $d$-th hyper-Kloosterman sum. As is well known, Kloosterman sums are trace functions  (see \cite{GKM}) and they satisfy Deligne's bound
$$|\Kl_d(n; q)|\leq d.$$ 
Therefore, possibly subject to the Ramanujan--Petersson conjecture, one obtains that \eqref{deltatransform} is bounded by
$$\frac{X}{q^{\frac{d+1}2}}\sum_{n}\lambda_{\ov\pi}(n)\Kl_d( a n;q)\widecheck{V}\left(\frac{nX}{q^d}\right)\ll  X^{o(1)}q^{\frac{d-1}{2}}$$
which improves over \eqref{deltasum} as long as
\begin{equation*}
	\label{arithprogrrange}
	q\leq X^{\theta_d-\delta},\ \delta>0
\end{equation*}
where \begin{equation}
	\label{thetaddef}
	\theta_d=\frac{2}{d+1}.
\end{equation}
We call the exponent $\theta_d$ the {\em standard level of distribution} of the sequence $(\lambda_\pi(n))_{n\geq 1}$. For instance we have
$$\theta_2=2/3,\ \theta_3=1/2,\ \theta_4=2/5,\ \theta_6=2/7.$$
To increase this standard level of distribution we would then need to obtain \eqref{nontrivialgeneric} for $K(n)=\Kl_d(an;q)$ and
$\widecheck{X}=q^d/X$ smaller than
$$\Xap:= q^{\frac{d-1}2}=\Xcv q^{-1/2}.$$
We call this the {\em a.p. range}.


Observe that the three results mentioned above  fall short of reaching the a.p. range.

One possible way to improve the situation is to exploit some special properties of the arithmetic function $\lambda_\pi$ like the existence of a Dirichlet factorisation
$$\lambda_\pi(n)=\lambda_{\pi_1}\star\lambda_{\pi_2}(n)=\sum_{l m=n}\lambda_{\pi_1}(l)\lambda_{\pi_2}(m).$$
A landmark example is the work of Friedlander and Iwaniec \cite{FI} on the ternary divisor function
$$d_3(n)=1\star 1\star 1(n)=\sum_{klm=n}1$$
where the standard exponent $\theta_3$ was replaced by $\theta_3+1/230$ (see \cite{HBActa,FKMMath,KMS2} for further improvements and generalisations of their ideas). See also \cite{FIActa} for result on higher order divisor functions. Recently, in joint work with E. Kowalski \cite{KLM}, we could pass the a.p. range $\theta_4=2/5$ for Rankin--Selberg coefficients $\lambda_{f\times f}(n)$ for $f$ a cusp form of level $1$ by using the factorisation
$$\lambda_{f\times f}(n)=1\star\lambda_{\sym^2f}(n)=\sum_{lm=n}\lambda_{\sym^2f}(m)$$
combined with \cite{KLMS}.

 \subsubsection*{Composite moduli}
Another option\footnote{Possibly combined with the previous one.} is to exploit existing {\em  factorisations} of the {\em modulus} $q$. This is for instance the case of the work of Fouvry, Iwaniec, and Katz \cite{FIK} on the divisor function $\lambda_\pi(n)=d_2(n)$ (which requires an additional averaging over one factor of the modulus  but extends to Fourier coefficients of cusp forms) or of Irving \cite{Irving} for sufficiently smooth moduli $q$ (however the method also uses the factorisation of the divisor function).

In this paper we improve the ranges \eqref{shortestd2} and \eqref{shortestd6} for suitable factorable moduli $q$. 
 
To simplify we limit ourselves to the case where the modulus is of the form $q=q_0q_1$  where $q_0$ and $q_1$ are distinct primes; consequently we assume that $$K:\Zz/q_0q_1\Zz\mapsto \Cc$$ can be expressed (via the Chinese Remainder Theorem) as a product of two functions of respective moduli $q_0$ and $q_1$.  For $K(\cdot)$ on $\Zz/q\Zz$ we denote
\begin{equation}\label{FourierTransform-K}
\what K(n)=\frac{1}{q^{1/2}}\sum_{x\in\Fq}K(x)e\left(\frac{nx}q\right)
\end{equation}
its normalized Fourier transform. We write $\|\widehat{K}\|_{\infty}$ for the
maximum of $|\what K(n)|$ for $n\in \Zz/q\Zz$.
 
Our first result is for $d=2$:
 \begin{theorem}\label{d2thm} Let $f(z)$ be a  Hecke eigencuspform  of level $1$ either holomorphic of weight $k\geq 2$ or a Laplacian eigenform with spectral parameter $t_f$; let $(\lambda_f(n))_{n\geq 1}$ be its Hecke eigenvalues. Let $q=q_0q_1$ be a product of two distinct primes and let $$K_0:\Zz/{q_0}\Zz\mapsto \Cc,\ K_1:\Zz/{q_1}\Zz\mapsto \Cc$$ be two complex valued functions which we identify with functions on $\Zz$ of period $q_0$ and $q_1$ respectively;  we assume that $K_0$ is the trace function attached to an $\ell$-adic middle extension Fourier sheaf $\mcF$ on $\Aa^1_{\Ff_{q_0}}$, geometrically irreducible, pure of weight $0$ and such that
 \begin{itemize}
 \item The automorphism group {\rm(}see \cite[\S 7.1]{FKMS}{\rm)} of the Fourier transform sheaf $\widehat\mcF$  is trivial.
     \item {\rm (MO)}  {\rm(}see \cite[\S 1]{LMS}{\rm)} There is no $ \lambda \in \Fqt$ such that the geometric monodromy group of $\mcF$ has some quotient which is equal, as a representation of the geometric fundamental group $\pi_1$ into an algebraic group, to the geometric monodromy group of the Kloosterman sheaf $[\times \lambda]^*\KL_2$ modulo $\pm 1$.

 \end{itemize}
 
 Let $K$ be the $q$-periodic function given by
$$
K(\cdot)=K_0(\cdot)K_1(\cdot).$$

Let $Z\geq 1$ be some parameter and $V\in \mcC^\infty_c(\Rr)$ be a smooth function compactly supported in the interval $[1,2[$ satisfying for all $j\geq 0$,
\begin{equation}\label{bound-of-V}
    V^{(j)}(x)\ll_j Z^j.
\end{equation}
For $X\geq 1$ we have the bound   
\begin{equation}
	\label{d2boundcomposite}
	\sum_{n=1}^{\infty}\lambda_f(n)K(n)V\left(\frac{n}{X}\right)\ll  \|\widehat{K_1}\|_{\infty} X^{o(1)}\left(  Z^{1/2}X^{1/2} q_0^{1/2}+ZX^{1/2}q^{1/2}{q_0}^{-1/4}+  Zq^{1/2}{q_0}^{1/4}\right),
\end{equation}
where the implicit constant depends at most on $f$, on the conductor $C_0=C(\mcF)$ of $\mcF$ and on the implicit constants in \eqref{bound-of-V}.
 \end{theorem}

\begin{remark} Examples of traces functions $K_0$ whose associated sheaf $\mcF$ satisfy the conditions of Theorem \ref{d2thm} include the hyper-Kloosterman sums $\Kl_d(\cdot;q_0)$ for $d\not=2$.
    
\end{remark}

 \begin{remark} 
 In particular, under the above assumptions on $K$, for $Z=1$ and
$$q_0= q^{2/3+o(1)}$$ 
from \eqref{d2boundcomposite} we see that the bound \eqref{nontrivialgeneric} holds as long as
\begin{equation*}\label{Weylrange}
	X\geq q^{2/3+\delta}
\end{equation*} for any given $\delta>0$. The range (for $d=2$) $$q^{2/3}=\Xcv^{2/3}>\Xap=q^{1/2}$$  is sometimes called the {\em Weyl range}. It can be related to the work of Heath-Brown  \cite{H-B78} who obtained a Weyl-type subconvexity bound $L(\chi,1/2)=O\big(q^{1/6+o(1)}\big)$ for the Dirichlet $L$-functions provided that $q$ has a factor $q_0$ of size $q_0\approx q^{2/3}$.

 \end{remark}
 \begin{remark}\label{van der Corput}
 It is possible to improve this range further assuming that $q_0$ is squarefree and suitably factorable (under additional assumptions on $K_0$) by using the
    recent work of Wu, Xi, and Sawin \cite{PW} which extends Heath-Brown's $q$-van der Corput's method to general trace functions.
\end{remark}
 
\begin{remark}
	In \cite{FKM1} the bound \eqref{nontrivialgeneric} was obtained for $q$ prime in the range \eqref{shortestd2} under the sole assumption that $\mcF$ is a Fourier sheaf (not geometrically isomorphic to the constant sheaf or any Artin--Schreier sheaf). Theorem \ref{d2thm} holds in this generality as well, but to simplify the exposition, we have chosen to make this extra assumption regarding the automorphism group of 
	$\widehat\mcF$.
\end{remark}


Our second main result concerns the case for $d=6$ where $\lambda_\pi(n)$ is given by \eqref{RS23def}, that is, the $n$-th coefficient of a $\GL_{3,\Qq}\times\GL_{2,\Qq}$ Rankin--Selberg $L$-function $L(\varphi\times f,s)$.

\begin{theorem}\label{mainthm2}  Let $f(z)$ be a  Hecke eigencuspform  of level $1$ either holomorphic of weight $k\geq 2$ or a Laplacian eigenform with spectral parameter $t_f$; let $(\lambda_f(n))_{n\geq 1}$ be its Hecke eigenvalues. Let $\varphi$ be a Hecke--Maass cuspform for $\GL_{3,\Qq}$ of level $1$ with Fourier coefficients $(\lambda_\varphi(n,r))_{n\geq 1,r\not=0}$.

Let $q=q_0q_1$ be a product of two distinct primes and let $$K_0:\Zz/{q_0}\Zz\mapsto \Cc,\ K_1:\Zz/{q_1}\Zz\mapsto \Cc$$ be two complex valued functions which we identify with functions on $\Zz$ of period $q_0$ and $q_1$ respectively;    we assume moreover that $K_0$ is the trace function attached to an $\ell$-adic middle extension Fourier sheaf $\mcF$ on $\Aa^1_{\Ff_{q_0}}$, geometrically irreducible and pure of weight $0$, of conductor $C_0$ and which is {\em good} in the sense of \cite[\S 1]{LMS}.  Let $K$ be the $q$-periodic function given by
$$
K(\cdot)=K_0(\cdot)K_1(\cdot).$$

Let $V\in \mcC^\infty_c(\Rr)$ be a smooth function satisfying \eqref{bound-of-V}. Let $X\geq 1$ be such that $X\geq Z^4q^2q_0^{1/2}$. Then
\begin{multline}
	\label{d3x2boundcomposite}
	\sum_{n\geq 1,r\not=0}\lambda_\varphi(n,r)\lambda_f(n)K(nr^2)V\left(\frac{nr^2}{X}\right)\ll_{\varphi,f,\|K_1\|_\infty,C_0}X^{o(1)}Z^2\Bigg(X^{3/4}{q_0}^{3/4}\\
+
X^{\frac{2-\theta_3}{3-2\theta_3}}(q^2q_0^{1/2})^{\frac{1-\theta_3}{3-2\theta_3}}+\frac{X}{{q_0}^{1/4}}+ \frac{X^{3/4}q}{q_0^{1/2}}\Bigg).
	\end{multline}
Here $\theta_3=5/14$ is the best known bound towards the Ramanujan--Petersson conjecture on $\GL_3$. 
\end{theorem}

 \begin{remark} In particular for $$q_0= q^{4/5+o(1)}$$ we obtain that
$$\sum_{n\geq 1,r\not=0}\lambda_\varphi(n,r)\lambda_f(n)K(nr^2)V\left(\frac{nr^2}{X}\right)\ll_{\varphi,f,\|K_1\|_\infty,C_0}X^{o(1)}Z^2\Bigg(X^{3/4}q^{3/5}+X^{\frac{2-\theta_3}{3-2\theta_3}}q^{\frac{12(1-\theta_3)}{5(3-2\theta_3)}}+\frac{X}{q^{1/5}}\Bigg)$$
and when $Z=1$ this bound is non-trivial {\rm(}i.e., \eqref{nontrivialgeneric} holds{\rm)} as soon as
$$X\geq q^{12/5+\delta}=q^{5/2-1/10+\delta}$$
for some $\delta>0$, in which case the second factor inside the parentheses can be removed.
\end{remark}

Applying Theorem \ref{mainthm2} to the function
$$K(n)=\Kl_6(an;q)$$
we obtain

\begin{corollary}\label{mainarith} 

Notations be as above. Assume that $$q_0= q^{4/5+o(1)}.$$
We have for any $(a,q)=1$,
$$\sum_\stacksum{n,|r|\geq 1}{nr^2\equiv a\bmod q}\lambda_\varphi(n,r)\lambda_f(n)V\left(\frac{nr^2}{X}\right)\ll_{\varphi,f} X^{o(1)}\left(X^{1/4}q^{8/5}+q^{23/10}\right).$$
In particular, given $\eta>0$,  we have
$$\sum_\stacksum{n,|r|\geq 1}{nr^2\equiv a\bmod q}\lambda_\varphi(n,r)\lambda_f(n)V\left(\frac{nr^2}{X}\right)\ll_{\varphi,f,\eta}  (\frac{X}{q})^{1-\delta}$$
for some $\delta=\delta(\eta)>0$ as long as
$$q\leq X^{2/7+1/364-\eta}=X^{\theta_6+1/364-\eta}.$$
\end{corollary}

\begin{remark}
The exponent $1/364$ matches the one presented in \cite[Cor. 1.4]{Lin-Sun}. In that work, Q. Sun and the first named author obtained the following error term for the ``sharp cut" sum of the coefficients $\sum_{mr^2=n}\lambda_\vphi(m,r)\lambda_f(m)$
$$\sum_{nr^2\leq X}\lambda_\varphi(n,r)\lambda_f(n)=O\left(X^{\tau_6-1/364+o(1)}\right)$$
where
$$\tau_d=\frac{d-1}{d+1}$$ is the ``standard" exponent of Landau (see \cite[Prop. 1.1]{FICan}) for sharp-cut sums of coefficients of automorphic $L$-functions of degree $d$.	
\end{remark}

\section{Proof of Theorem \ref{d2thm}}
 The proof of Theorem \ref{d2thm} is based on the method introduced in \cite{AHLS}, in which a Burgess type bound was obtained for $K=\chi$ a multiplicative character of prime modulus $q$, by unexpectedly utilizing a {\em ``trivial" delta method}.  We demonstrate in this work once again that the trivial delta method is strong enough to reach a Weyl type bound when the moduli $q$ factors in an appropriate way (see \cites{subGL2GL2,Agg1} for other interesting applications of the trivial delta symbol).
 
 Let $(\lambda_f(n))_{n\geq 1}$ be the Hecke eigenvalues of a $\rm GL_2$ cusp form $f$. Let $q={q_0}{q_1}$. Let $K \bmod q$ be a trace function. 
 From the definition \eqref{FourierTransform-K} we have the twisted multiplicativity 
\begin{equation}\label{twisted-multiplicity}
    \widehat{K}(b)=\widehat{K_0}(\overline{{q_1}}b)\widehat{K_1}(\overline{{q_0}}b).
\end{equation}
 Consider 
\begin{equation}
    S(X):=\sum_{n=1}^{\infty}\lambda_f(n)K(n)V\left(\frac{n}{X}\right),
\end{equation}
where $V$ satisfies the bound \eqref{bound-of-V}.
Then
\begin{equation*}
    S(X)=\sum_{n=1}^{\infty}\lambda_f(n)W\left(\frac{n}{X}\right)\sum_{r=1}^{\infty}K(r)V\left(\frac{r}{X}\right)\delta_{n=r}.
\end{equation*}
Here $W$ is a smooth function supported in $(1/100,100)$ and satisfying $W^{(j)}(x)\ll_j 1$ for $j\geq 0$ and $W(x)=1$ for $x\in [1,2]$, and $\delta_{n=r}$ is the Kronecker delta symbol.

Following \cite{AHLS}, if we assume that $p$ is a prime coprime with $q$ and such that
\begin{equation}\label{initial-assumption}
    p{q_0}>1000X
\end{equation} then for $|n-r|\leq 100X$, we have
\begin{equation*}
    \begin{split}
   \delta_{n=r}=&\frac{1}{p{q_0}}\sum_{u(p{q_0})}
   e\left(\frac{u(n-r)}{p{q_0}}\right)\\
   =&\frac{1}{p{q_0}}\sum_{c|p{q_0}}\sumstar_{\alpha(c)}e\left(\frac{\alpha(n-r)}{c}\right),
\end{split}\end{equation*}
and by summing over all $p\in {\rm P}:=\{p\in [P,2P[: p\, \text{prime},\ (p,q)=1\}$, we have
\begin{equation*}
   \delta_{n=r}=\frac{1}{|\rm P|}\sum_{p\in \rm P}\frac{1}{pq_0}\sum_{c|pq_0}\sumstar_{\alpha(c)}e\left(\frac{\alpha(n-r)}{c}\right).
\end{equation*}
with $|{\rm P}|\gg P/\log P$ from the prime number theorem.
 \begin{remark} Choosing the parameter $q_0$ in place of $q$ is one chief difference between this paper and \cite{AHLS}.
\end{remark}
Applying the above expression of  
$\delta_{n=r}$ to the sum $S(X)$, we get
\begin{equation*}
    S(X)=\frac{1}{|{\rm P}|q_0}\sum_{p\in \rm P}\frac{1}{p}\sum_{c|pq_0}\sumstar_{\alpha(c)}\sum_{n=1}^{\infty}\lambda_f(n)e\left(\frac{\alpha n}{c}\right)W\left(\frac{n}{X}\right)\sum_{r=1}^{\infty}K(r)e\left(\frac{-\alpha r}{c}\right)V\left(\frac{r}{X}\right).
\end{equation*}
Applying the Voronoi summation formula (\cite[Thm. A.4]{KMV}), we obtain
\begin{equation*}
\sum_{n=1}^{\infty}\lambda_f(n)e\left(\frac{\alpha n}{c}\right)W\left(\frac{n}{X}\right)=\frac{X}{c}\sum_{\pm}\sum_{n=1}^{\infty}\overline{\lambda_f(n)}e\left(\frac{\mp\bar{\alpha} n}{c}\right)\widehat{W}^{\pm}\left(\frac{n}{c^2/X}\right).
\end{equation*}
Here $\widehat{W}^{\pm}(\cdot)$ are defined as in \cite[Lem. 2.2]{AHLS}.
Applying the Poisson summation formula (modulo $[c,q]$), the $r$-sum
\begin{equation*}
\sum_{r\geq1}K(r)e\left(\frac{-\alpha r}{c}\right)V\left(\frac{r}{X}\right)
\end{equation*}
is transformed into
\begin{equation}\label{r-sum}
\frac{X}{[c,q]} \sum_{r\in\mathbb{Z}}\left(\sum_{\beta\bmod[c,q]}K(\beta)e\left(\frac{-\alpha\beta}{c}\right)e\left(\frac{r\beta}{[c,q]}\right)\right)\widehat{V}\left(\frac{rX}{[c,q]}\right),
\end{equation}
where 
$\widehat{V}$ denotes the Fourier transform of the function $V$.

Following \cite{AHLS}, we introduce the useful notation $a_b:=\frac{a}{(a,b)}$. Then, using the relation $c=\frac{qc_q}{q_c}$ we can write $$e\left(\frac{\bullet}{c}\right)=e\left(\frac{\bullet q_c}{qc_q}\right)=e\left(\frac{\bullet q_c\overline{c_q}}{q}\right)e\left(\frac{\bullet q_c\bar{q}}{c_q}\right).$$
Similarly, using the relation $[c,q]=qc_q$ we have the reciprocity relation
 $$e\left(\frac{r\beta}{[c,q]}\right)=e\left(\frac{r\beta\overline{c_q}}{q}\right)e\left(\frac{r\beta\bar{q}}{c_q}\right).$$ 
Then the $\beta$-sum in \eqref{r-sum} can be rewritten as
\begin{equation*}
\begin{split}
&\sum_{\beta\bmod q}K(\beta)e\left(\frac{(r-\alpha q_c)\overline{c_q}\beta}{q}\right) \times \sum_{\beta\bmod c_q}e\left(\frac{(r-\alpha q_c)\overline{q}\beta}{c_q}\right)\\
&= q^{1/2}c_q\widehat{K}((r-\alpha q_c)\overline{c_q})\times \delta_{r-\alpha q_c\equiv0\bmod c_q},
\end{split}
\end{equation*}
where $\widehat{K}$ is the normalized Fourier transform \eqref{FourierTransform-K} and $\overline{c_q}$ is the multiplicative inverse of $c_q$ modulo $q$. We also note that the weight function $\widehat{V}\left(\frac{rX}{[c,q]}\right)$ restricts the effective range of the $r$-sum in \eqref{r-sum} to 
$$|r|\leq X^{o(1)}Z[c,q]/X.$$

 By this we mean that for any $\eps>0$, the contribution of the terms satisfying $$|r|> X^{\eps}Z[c,q]/X$$
is bounded by $O_{A,\eps}(X^{-A})$ for any $A\geq 1$.

\begin{notation} To lighten the expressions to come we write 
$$A\lesssim B\hbox{ in place of }A\leq X^{o(1)}B.$$
\end{notation}

Therefore \eqref{r-sum} becomes
\begin{equation*}
\frac{X}{q^{1/2}}\underset{\begin{subarray}{c}
|r|\lesssim Z[c,q]/X \\ r-\alpha q_c\equiv0\bmod c_q
\end{subarray}}{\sum}\widehat{K}((r-\alpha q_c)\overline{c_q}) \widehat{V}\left(\frac{rX}{[c,q]}\right) + O_A\left(X^{-A}\right).
\end{equation*}

Therefore one obtains that
\begin{multline*}
    S(X)=\frac{1}{|{\rm P}|q_0}\sum_{p\in \rm P}\frac{1}{p}\sum_{c|pq_0}\sumstar_{\alpha(c)}\frac{X}{c}\sum_{\pm}\sum_{n=1}^{\infty}\overline{\lambda_f(n)}e\left(\frac{\mp\bar{\alpha} n}{c}\right)\widehat{W}^{\pm}\left(\frac{n}{c^2/X}\right)\\
\times\frac{X}{q^{1/2}}\underset{\begin{subarray}{c}
|r|\lesssim Z[c,q]/X \\ r-\alpha q_c\equiv0\bmod c_q
\end{subarray}}{\sum}\widehat{K}((r-\alpha q_c)\overline{c_q}) \widehat{V}\left(\frac{rX}{[c,q]}\right)+ O_A\left(X^{-A}\right).
\end{multline*}

If $c=pq_0$, then in the above notation $c_q=p$ and $q_c=q_1$. The congruence condition $r-\alpha q_c\equiv0\bmod c_q$ becomes $\bar{\alpha}\equiv \bar{r}q_1\bmod p$.
Furthermore, from the twisted multiplicativity in \eqref{twisted-multiplicity}, we have
\begin{equation*}
\begin{split}
\widehat{K}((r-\alpha {q_1})\bar{p})
=\widehat{K_0}((r-\alpha {q_1})\overline{{q_1}p})\widehat{K_1}(r\overline{{q_0}p});
\end{split}\end{equation*}
similarly we have 
\begin{equation*}
\begin{split}
e\left(\frac{\mp\bar{\alpha} n}{pq_0}\right)=e\left(\frac{\mp\bar{\alpha} n\overline{q_0}}{p}\right)e\left(\frac{\mp\bar{\alpha} n\bar{p}}{q_0}\right).
\end{split}\end{equation*}

Therefore we can further rewrite $S(X)$ as
\begin{multline}\label{main+errors}
   S(X)=\frac{X^2}{|{\rm P}|{q_0}^2q^{1/2} }\sum_{\pm}\sum_{n=1}^{\infty}\overline{\lambda_f(n)}\sum_{p\in \rm P}\frac{1}{p^2}\sum_{|r|\lesssim  Zpq/X}\widehat{K_1}(r\overline{{q_0}p})e\left(\frac{\mp\overline{r{q_0}} n{q_1}}{p}\right)\widehat{V}\left(\frac{rX}{pq}\right)\\
\times \sumstar_{\alpha({q_0})}\widehat{K_0}((r-\alpha {q_1})\overline{{q_1}p})e\left(\frac{\mp\bar{\alpha} n\bar{p}}{{q_0}}\right)\widehat{W}^{\pm}\left(\frac{n}{p^2{q_0}^2/X}\right) +S(X)_{c=p}+S(X)_{c=q_0}+ O\left(X^{-A}\right),
\end{multline}
where the term $S(X)_{c=p},S(X)_{c=q_0}$, given in \eqref{c=p-term} and \eqref{c=q0-term}, correspond to the contribution from the terms $c=p$
and $c=q_0$ respectively (we notice that due to rapid decay of $\widehat{W}^{\pm}(x)$ the $c=1$ term was absorbed into the negligible error term $O_A\left(X^{-A}\right)$). We will treat such terms separately later in Sec. \ref{errorterm-treatment} (see the bounds \eqref{error-c=p} and \eqref{error-c=q0}).

\begin{remark}
Estimating trivially (assuming the sum over $\alpha\bmod q_0$ in \eqref{main+errors} is 
 bounded above by $O(q_0^{1/2})$) and ignoring the error terms, we have
\begin{equation*}
\begin{split}
    S(X)\ll& \frac{X^{2+o(1)}}{P^3{q_0}^2q^{1/2} }\frac{p^2{q_0}^2}{X}P\frac{Zpq}{X}q_0^{1/2}
    \ll X^{o(1)}ZPq^{1/2}q_0^{1/2}.
    \end{split}\end{equation*}
    In view of the constraint \eqref{initial-assumption}, this bound is insufficient to beat the trivial bound $O(X)$ yet.
\end{remark}
We can effectively truncate the $n$ sum at $n\lesssim P^2{q_0}^2/X$. To this end, we further break the $n$-sum into dyadic intervals by introducing another smooth weight $U(\frac{n}{N})$ to the $n$-sum, where the $N$'s are such that
$$N\lesssim \frac{P^2{q_0}^2}{X}.$$
Moreover we can pull out the factor $\frac{1}{p^2}$ from the $p$-sum by introducing another smooth weight to the $p$-sum which we suppress from our notation. To conclude, we can bound $S(X)$ as follows
\begin{multline*}
S(X)\ll \frac{X^{2+o(1)}}{P^3{q_0}^2q^{1/2} }\sup_{N\lesssim P^2{q_0}^2/X}\bigg|\sum_{n=1}^{\infty}\overline{\lambda_f(n)}U\left(\frac{n}{N}\right)\sum_{p\in \rm P}\sum_{|r|\lesssim Zpq/X}\widehat{K_1}(r\overline{{q_0}p})e\left(\frac{\mp\overline{r{q_0}} n{q_1}}{p}\right)\widehat{V}\left(\frac{rX}{pq}\right)\\
\times\sumstar_{\alpha({q_0})}\widehat{K_0}((r-\alpha {q_1})\overline{{q_1}p})e\left(\frac{\mp\bar{\alpha} n\bar{p}}{{q_0}}\right) \widehat{W}^{\pm}\left(\frac{n}{p^2{q_0}^2/X}\right)\bigg|+|S(X)_{c=p}|+|S(X)_{c=q_0}|.
\end{multline*}

Applying Cauchy--Schwarz inequality, one has
\begin{multline*}
    S(X)\ll |S(X)_{c=p}|+|S(X)_{c=q_0}|+\frac{X^{2+o(1)}}{P^3{q_0}^2q^{1/2} }\sup_{N\lesssim P^2{q_0}^2/X}\bigg(\sum_{n=1}^{\infty}|\overline{\lambda_f(n)}|^2U\left(\frac{n}{N}\right)\bigg)^{1/2}\\
    \times\bigg(\sum_{n=1}^{\infty}\bigg|\sum_{p\in \rm P}\sum_{|r|\lesssim Zpq/X}\widehat{K_1}(r\overline{{q_0}p})e\left(\frac{\mp\overline{r{q_0}} n{q_1}}{p}\right)\widehat{V}\left(\frac{rX}{pq}\right)\\
    \times\sumstar_{\alpha({q_0})}\widehat{K_0}((r-\alpha {q_1})\overline{{q_1}p})e\left(\frac{\mp\bar{\alpha} n\bar{p}}{{q_0}}\right)\widehat{W}^{\pm}\left(\frac{n}{p^2{q_0}^2/X}\right)\bigg|^2U\left(\frac{n}{N}\right)\bigg)^{1/2}.
  \end{multline*}
Now we proceed to open the square and obtain that
\begin{align*}
	 S(X)\ll& |S(X)_{c=p}|+|S(X)_{c=q_0}|+\frac{X^{2+o(1)}}{P^3{q_0}^2q^{1/2} }\sup_{N\lesssim P^2{q_0}^2/X}N^{1/2}\bigg(\sum_{p_1,p_2\in \rmP}\\
	 &\times \sum_{|r_1|,|r_2|\lesssim ZPq/X}\widehat{K_1}(r_1\overline{{q_0}p_1})\overline{\widehat{K_1}(r_2\overline{{q_0}p_2})}\widehat{V}\left(\frac{r_1X}{p_1q}\right)\overline{\widehat{V}\left(\frac{r_2X}{p_2q}\right)}\\
	 &\times\sumstar_{\alpha_1({q_0})}\sumstar_{\alpha_2({q_0})}\widehat{K_0}((r_1-\alpha_1 {q_1})\overline{{q_1}p_1})\overline{\widehat{K_0}((r_2-\alpha_2 {q_1})\overline{{q_1}p_2})}\\
   &\times \sum_{n=1}^{\infty}e\left(\mp\frac{\overline{r_1{q_0}} n{q_1}p_2-\overline{r_2{q_0}} n{q_1}p_1}{p_1p_2}\right)e\left(\mp\frac{\overline{\alpha_1} n\overline{p_1}-\overline{\alpha_2} n\overline{p_2}}{{q_0}}\right) \mcW\left(\frac{n}{N}\right)\bigg)^{1/2}.
\end{align*}
where 
$$\mcW\left(\frac{n}{N}\right):=U\left(\frac{n}{N}\right)\widehat{W}^{\pm}\left(\frac{n}{{p^2_1}{q_0}^2/X}\right)\overline{\widehat{W}^{\pm}}\left(\frac{n}{{p^2_2}{q_0}^2/X}\right).$$
Applying the Poisson summation formula, one has  
\begin{multline*}
\sum_{n=1}^{\infty}e\left(\mp\frac{\overline{r_1{q_0}} n{q_1}p_2-\overline{r_2{q_0}} n{q_1}p_1}{p_1p_2}\right)e\left(\mp\frac{\overline{\alpha_1} n\overline{p_1}-\overline{\alpha_2} n\overline{p_2}}{{q_0}}\right) \mcW\left(\frac{n}{N}\right)\\
=\frac{N}{p_1p_2{q_0}}\sum_{\tilde{n}\in \Zz}\sum_{\beta(p_1p_2)}e\left(\mp\frac{(\overline{r_1{q_0}} {q_1}p_2-\overline{r_2{q_0}} {q_1}p_1)\beta\mp\tilde{n}\beta\overline{{q_0}}}{p_1p_2}\right)\\
\times\sum_{\beta({q_0})}e\left(\mp\frac{(\overline{\alpha_1p_1}-\overline{\alpha_2p_2})\beta\mp\tilde{n}\beta\overline{p_1p_2}}{{q_0}}\right) \widehat{\mcW}\left(\frac{\tilde{n}N}{p_1p_2{q_0}}\right)\\
=N\sum_{\tilde{n}\in \Zz}\delta_{\alpha_2\equiv \overline{\overline{\alpha_1}p_2\mp\tilde{n}}\,p_1\bmod  {q_0}}\delta_{\tilde{n}\equiv \mp(\overline{r_1}p_2-\overline{r_2} p_1){q_1}\bmod  {p_1p_2}}\widehat{\mcW}\left(\frac{\tilde{n}N}{p_1p_2{q_0}}\right).
    \end{multline*}
Here $\widehat{\mcW}$ denotes the Fourier transform of $\mcW$. Plugging these calculations back, we find
\begin{multline*}
    S(X)\ll |S(X)_{c=p}|+|S(X)_{c=q_0}|+\frac{X^{2+o(1)}}{P^3{q_0}^2q^{1/2} }\sup_{N\lesssim P^2{q_0}^2/X}N\bigg(\sum_{p_1,p_2\in \rmP}\sum_{|r_1|,|r_2|\lesssim ZPq/X}\widehat{K_1}(r_1\overline{{q_0}p_1})\\
    \times \overline{\widehat{K_1}(r_2\overline{{q_0}p_2})}\widehat{V}\left(\frac{r_1X}{p_1q}\right)\overline{\widehat{V}\left(\frac{r_2X}{p_2q}\right)}q_0^{1/2}\mathop{\sum_{\tilde{n}\in \Zz}}_{\tilde{n}\equiv \mp(\overline{r_1}p_2-\overline{r_2} p_1){q_1}\bmod p_1p_2} \mathfrak{C}_{\gamma,{q_0}}(\tilde{n})\widehat{\mcW}\left(\frac{\tilde{n}N}{p_1p_2{q_0}}\right)\bigg)^{1/2}
    \end{multline*}
    with $\mathfrak{C}_{\gamma,{q_0}}(\tilde{n})$ the following correlation sum
    \begin{equation}\label{correlation-sum}
\begin{split}
   \mathfrak{C}_{\gamma,{q_0}}(\tilde{n}):=&\frac{1}{q_0^{1/2}}\sumstar_{\alpha({q_0})}\widehat{K_0}((r_1-\alpha {q_1})\overline{{q_1}p_1})
   \overline{\widehat{K_0}((r_2- \overline{\bar{\alpha}p_2\mp\tilde{n}}\,p_1 {q_1})\overline{{q_1}p_2})}\\
   =&\frac{1}{q_0^{1/2}}\sumstar_{\alpha({q_0})}\widehat{K_0}((\overline{{q_1}}r_1-\alpha )\overline{p_1})
   \overline{\widehat{K_0}((\overline{{q_1}}r_2- \overline{\bar{\alpha}p_2\mp\tilde{n}}\,p_1 )\overline{p_2}))}.
    \end{split}\end{equation}
     \begin{remark}
As can be expected, if the modulus $q=q_0q_1$ of $K(\cdot)$ admits further factorisation, then it is possible to apply a $q$-analogue of van der Corput's argument to the $(r_1,r_2)$-sum to obtain further improvements of our final result, cf. see also Remark \ref{van der Corput}.
\end{remark}
  We rewrite the inequality above into the form  
    \begin{equation}\label{final-object}
    S(X)\ll 
    |S(X)_{c=p}|+|S(X)_{c=q_0}|+S_{\text{dia}}(X)+S_{\text{off}}(X),
    \end{equation}
    where
    $S_{\text{dia}}(X)$ and $S_{\text{off}}(X)$ denote the terms corresponding to $\tilde{n}=0$ and $\tilde{n}\neq 0$ respectively.

\subsection{Contribution from the $\tilde{n}\not=0$ frequencies}
  Now we treat the term $S_{\text{off}}(X)$, corresponding to $\tilde{n}\neq 0$ in \eqref{final-object}.
  
   We consider two subcases:  $\tilde{n}\not\equiv 0\bmod q_0$ and $\tilde{n}\equiv 0\bmod q_0$ and denote their contribution to $S_{\text{off}}(X)$ by $S_{\text{off},1}(X)$ and $S_{\text{off},2}(X)$, respectively.

   \subsubsection{The case  $\tilde{n}\not\equiv 0\bmod q_0$}
   By a change of variable, the correlation sum
    can be written
   $$\mathfrak{C}_{\gamma,{q_0}}(\tilde{n})=\frac{1}{q_0^{1/2}}\sumstar_{\alpha({q_0})}\widehat{K_0}(\gamma_1\cdot\alpha)
   \overline{\widehat{K_0}(\gamma_2\cdot\alpha)}=\frac{1}{q_0^{1/2}}\sumstar_{\alpha({q_0})}\widehat{K_0}(\alpha)
   \overline{\widehat{K_0}(\gamma_2.\gamma_1^{-1}\cdot\alpha)}+O(\frac{\|\what K_0\|_\infty^2}{q_0^{1/2}})$$
   where
   $$\gamma_1=
   \begin{pmatrix}-q_1&r_1\\0&p_1q_1
   \end{pmatrix},\ \gamma_2=
   \begin{pmatrix}\mp \tilde n r_2-p_1q_1&p_2r_2\\
   \mp\tilde n p_2q_1&p_2^2q_1
   \end{pmatrix}\in\GL_2(\Ff_{q_0})
   $$
and where we have set for
$$   \gamma=\begin{pmatrix}
       a&b\\c&d
   \end{pmatrix}\in\GL_2(\Ff_{q_0}),\ \gamma\cdot\alpha=\frac{a\alpha+b}{c\alpha+d}.$$
If $\tilde n\not\equiv 0\mods{q_0}$, the 
matrix $\gamma_2.\gamma_1^{-1}$ is not a scalar matrix. Therefore, by our assumption that the group of automorphisms of the Fourier transform sheaf $\what\mcF$ is trivial and by \cite[Prop. 7.2]{FKMS}, we obtain 
        \begin{equation}\label{condition2-on-K0} \mathfrak{C}_{\gamma,{q_0}}(\tilde{n})\ll_{C_0} 1 \quad \text{whenever}\,\, \tilde{n}\not\equiv 0\bmod q_0.
    \end{equation}
  It follows that
   \begin{multline*}
       S_{\text{off},1}(X)\ll \frac{X^{2+o(1)}}{P^3{q_0}^2q^{1/2} }\times\\
       \sup_{N\lesssim  \frac{P^2{q_0}^2}X}N\bigg(q_0^{1/2}\sum_{p_1,p_2\in \rmP}\sum_{0\neq |\tilde{n}|\lesssim P^2{q_0}/N}\sum_\stacksum{|r_1|\lesssim ZPq/X}{r_1\equiv \mp\ov{\tilde{n}}p_2q_1\bmod p_1}\sum_\stacksum{|r_2|\lesssim ZPq/X}{r_2\equiv \pm\ov{\tilde{n}}p_1q_1\bmod p_2}\|\widehat{K_1}\|_{\infty}^2
    \bigg)^{1/2}
   \end{multline*}
which gives
   \begin{equation*}
\begin{split}
   S_{\text{off},1}(X) \ll&\frac{X^{2+o(1)}}{P^3{q_0}^2q^{1/2} }\sup_{N\lesssim P^2{q_0}^2/X}N\bigg(q_0^{1/2}P^2\frac{P^2{q_0}}{N}\frac{Zq}{X}\frac{Zq}{X}\|\widehat{K_1}\|_{\infty}^2
    \bigg)^{1/2}\\
    \ll&\frac{\|\widehat{K_1}\|_{\infty}Z X^{1+o(1)}q^{1/2}}{P{q_0}^{5/4}}\sup_{N\lesssim  P^2{q_0}^2/X}N^{1/2}\\
    \ll&\|\widehat{K_1}\|_{\infty}Z X^{1/2+o(1)}q^{1/2}{q_0}^{-1/4}.
    \end{split}\end{equation*}
    
       \subsubsection{The case  $\tilde{n}\equiv 0\mods{q_0}$}
       
       If $\tilde{n}\equiv 0\mods{q_0}$ we have
       
    \begin{equation}\label{correlation-sum-0}
       \mathfrak{C}_{\gamma,{q_0}}(0)
   =\frac{1}{q_0^{1/2}}\sumstar_{\alpha({q_0})}\widehat{K_0}(\alpha)
   \overline{\widehat{K_0}(\begin{pmatrix}
       p_1^2q_1&-(p_1r_1-p_2r_2)\\0&p_2^2q_1
   \end{pmatrix}\cdot\alpha)}+O(\frac{\|\what K_0\|_\infty^2}{q_0^{1/2}}).
    \end{equation}
    The matrix $\begin{pmatrix}
       p_1^2q_1&-(p_1r_1-p_2r_2)\\0&p_2^2q_1
   \end{pmatrix}$ is scalar if and only if
   $$p_1r_1-p_2r_2\equiv 0\mods{q_0}\hbox{ and }p_1\equiv \pm p_2\mods{q_0}.$$
   Again, by our assumption that the group of automorphisms of the Fourier transform sheaf $\what\mcF$ is trivial and by \cite[Prop. 7.2]{FKMS}, we obtain
    \begin{equation}\label{bound-when-0}
    \mathfrak{C}_{\gamma,{q_0}}(0)\ll_{C_0} 
q_0^{1/2}\delta_{\substack{r_1p_1\equiv r_2p_2\mods{q_0}\\ p_1\equiv \pm p_2\mods{q_0}}}+1.
\end{equation}

       Hence for the terms such that $\tilde{n}\equiv 0\bmod q_0$, we replace the previous $\tilde{n}$ by $q_0\tilde{n}$ with $\tilde{n}\lesssim 1+P^2/N$ and apply the bound \eqref{bound-when-0} to obtain
          \begin{equation*}
\begin{split}
    S_{\text{off},2}(X)\ll& \frac{X^{2+o(1)}}{P^3{q_0}^2q^{1/2} }\sup_{N\lesssim \frac{P^2{q_0}^2}{X}}N\bigg(q_0^{1/2}\sum_{p_1,p_2\in \rmP}\sum_{0< |\tilde{n}|\lesssim 1+P^2/N}\\
   & \ \ \ \times\sum_\stacksum{|r_1|\lesssim ZPq/X}{r_1\equiv \mp\ov{q_0\tilde{n}}p_2q_1\bmod p_1}\sum_\stacksum{|r_2|\lesssim ZPq/X}{r_2\equiv \pm\ov{q_0\tilde{n}}p_1q_1\bmod p_2}\|\widehat{K_1}\|_{\infty}^2\big(q_0^{1/2}\delta_{r_1p_1\equiv r_2p_2\bmod  {q_0}}+1\big)
    \bigg)^{1/2}\\
    \ll&\|\widehat{K_1}\|_{\infty}\frac{X^{2+o(1)}}{P^3{q_0}^2q^{1/2} }\sup_{N\lesssim P^2{q_0}^2/X}N\bigg(q_0^{1/2}P^2\big(1+\frac{P^2}{N}\big)\frac{Zq}{X}\frac{Zq}{X}
    \bigg)^{1/2}\\
    \ll&\|\widehat{K_1}\|_{\infty} ZX^{o(1)}\Bigl(q^{1/2}{q_0}^{1/4}+X^{1/2}q^{1/2}{q_0}^{-3/4}\Bigr).
    \end{split}\end{equation*}
 Combining the above estimates for $S_{\text{off},1}(X)$ and $S_{\text{off},2}(X)$, we have therefore proved that the $\tilde{n}\not=0$ frequencies $S_{\text{dia}}(X)$ in \eqref{final-object} contribute at most
      \begin{equation}\label{S-off}
\begin{split}
S_{\text{off}}(X)
\ll X^{o(1)}\|\widehat{K_1}\|_{\infty} \Bigl(ZX^{1/2}q^{1/2}{q_0}^{-1/4}+ Zq^{1/2}{q_0}^{1/4}\Bigr).
    \end{split}\end{equation}

       \subsection{Contribution of the $\tilde{n}=0$ frequency}
       We recall that the $\tilde{n}=0$ frequency $S_{\text{dia}}(X)$ in \eqref{final-object} is given by
       \begin{multline*}
    S_{\text{dia}}(X)= \frac{X^{2+o(1)}}{P^3{q_0}^2q^{1/2} }\sup_{N\lesssim P^2{q_0}^2/X}N\bigg(\sum_{p_1,p_2\in \rmP}\sum_{|r_1|,|r_2|\lesssim ZPq/X}\widehat{K_1}(r_1\overline{{q_0}p_1})\overline{\widehat{K_1}(r_2\overline{{q_0}p_2})}\\
    \times \widehat{V}\left(\frac{r_1X}{p_1q}\right)\overline{\widehat{V}\left(\frac{r_2X}{p_2q}\right)}\,q_0^{1/2} \mathfrak{C}_{\gamma,{q_0}}(0)\delta_{\overline{r_1}p_2\equiv\overline{r_2} p_1\bmod  {p_1p_2}}\widehat{\mcW}\left(0\right)\bigg)^{1/2}    \end{multline*}
       
  The congruence condition $\overline{r_1}p_2\equiv\overline{r_2} p_1\bmod p_1p_2$ implies that    $$p_1=p_2:=p\hbox{ and }r_1\equiv r_2\bmod p.$$
  Inserting $p_1=p_2$ into \eqref{correlation-sum-0}, the bound in  \eqref{bound-when-0}  reads
$$ \mathfrak{C}_{\gamma,{q_0}}(0)\ll_{C_0} 
q_0^{1/2}\delta_{r_1\equiv r_2\mods{q_0}}+1.$$
It follows that
         \begin{equation}\label{S-dia}
\begin{split}
    S_{\text{dia}}(X)\ll_{C_0} & \frac{X^{2+o(1)}}{P^3{q_0}^{2}q^{1/2} }\sup_{N\lesssim P^2{q_0}^2/X}N\bigg(q_0^{1/2}\sum_{p\in \rmP}\sum_{|r_1|,|r_2|\lesssim  ZPq/X}\|\widehat{K_1}\|_{\infty}^2\\
    &\times \big(q_0^{1/2}\delta_{r_1\equiv r_2\bmod  {q_0}}+1\big)\delta_{r_1\equiv r_2\bmod  p}\bigg)^{1/2}\\
    \ll_{C_0}&  \|\widehat{K_1}\|_{\infty}\frac{X^{2+o(1)}}{P^3{q_0}^{7/4}q^{1/2} }\frac{P^2{q_0}^2}{X}\bigg(\sum_{p\in \rmP}\sum_{|r_1|,|r_2|\lesssim  ZPq/X}q_0^{1/2}\delta_{r_1\equiv r_2\bmod  {pq_0}}\\
    &+\sum_{p\in \rmP}\sum_{|r_1|,|r_2|\lesssim  ZPq/X}\delta_{r_1\equiv r_2\bmod  p}\bigg)^{1/2}\\
    \ll_{C_0}&    \|\widehat{K_1}\|_{\infty} \frac{X^{1+o(1)}{q_0}^{1/4}}{Pq^{1/2} }\bigg(\frac{Z^{1/2} Pq^{1/2}{q_0}^{1/4}}{X^{1/2}}+\frac{ZPq}{X}\bigg)\\
    \ll_{C_0}& \|\widehat{K_1}\|_{\infty}  Z^{1/2}X^{1/2+o(1)} q_0^{1/2}+ \|\widehat{K_1}\|_{\infty} X^{o(1)}Zq^{1/2}{q_0}^{1/4}.
    \end{split}\end{equation}     
       In the above we have applied the fact that when $r_1\equiv r_2\mods{pq_0}$ then $r_1=r_2$, as we have (or we will) assume
       $$Zq<Xq_0.$$
    
\subsection{Error terms: treatment of $S(X)_{c=p}$ and $S(X)_{c=q_0}$}\label{errorterm-treatment}
Recall that in \eqref{main+errors} $S(X)_{c=p}$ is defined by
\begin{equation}\label{c=p-term}
\begin{split}
S(X)_{c=p}=&\frac{X}{|{\rm P}|q_0}\sum_{\pm}\sum_{p\in \rm P}\frac{1}{p^2}\sumstar_{\alpha(p)}\sum_{n=1}^{\infty}\overline{\lambda_f(n)}e\left(\frac{\mp\bar{\alpha} n}{p}\right)\widehat{W}^{\pm}\left(\frac{n}{p^2/X}\right) \times\frac{X}{q^{1/2}}\underset{\begin{subarray}{c}
|r|\lesssim Zpq/X \\ r-\alpha q\equiv0\bmod p
\end{subarray}}{\sum}\widehat{K}(r\bar{p};q) \widehat{V}\left(\frac{rX}{pq}\right)\\
=&\frac{X^2}{|{\rm P}|q_0q^{1/2}}\sum_{\pm}\sum_{p\in \rm P}\frac{1}{p^2}\sum_{|r|\lesssim  Zpq/X}\widehat{K}(r\bar{p};q)\widehat{V}\left(\frac{rX}{pq}\right)\sum_{n=1}^{\infty}\overline{\lambda_f(n)}e\left(\frac{\mp q\bar{r} n}{p}\right)\widehat{W}^{\pm}\left(\frac{n}{p^2/X}\right).
\end{split}\end{equation}
If we apply summation by parts and the Wilton-type bound (\cite[Thm. 8.1]{Iwaniec-Spec}) to the $n$-sum, we get
\begin{equation}\label{error-c=p}
\begin{split}
|S(X)_{c=p}|
\ll&\frac{X^{2+o(1)}}{Pq_0q^{1/2}}\frac{1}{P} \frac{ZPq}{X}\|\widehat{K}\|_{\infty}\big(\frac{P^2}{X}\big)^{1/2}\ll \frac{\|\widehat{K}\|_{\infty} ZX^{1/2+o(1)}q^{1/2}}{q_0}.
\end{split}\end{equation}

Likewise, $S(X)_{c=q_0}$ is given by
\begin{equation}\label{c=q0-term}
\begin{split}
S(X)_{c=q_0}=&\frac{X}{|{\rm P}|{q_0}^2}\sum_{\pm}\sum_{p\in \rm P}\frac{1}{p}\sumstar_{\alpha(q_0)}\sum_{n=1}^{\infty}\overline{\lambda_f(n)}e\left(\frac{\mp\bar{\alpha} n}{q_0}\right)\widehat{W}^{\pm}\left(\frac{n}{{q_0}^2/X}\right)\\
    &\quad\quad\quad \times\frac{X}{q^{1/2}}\sum_{|r|\lesssim Zq/X}\widehat{K}(r-\alpha q_1;q) \widehat{V}\left(\frac{rX}{q}\right)\\
    =&\frac{X^2}{|{\rm P}|{q_0}^{3/2}q^{1/2}}\sum_{\pm}\sum_{p\in \rm P}\frac{1}{p}\sum_{n=1}^{\infty}\overline{\lambda_f(n)}\widehat{W}^{\pm}\left(\frac{n}{{q_0}^2/X}\right)\\
    &\quad\quad\quad \times\sum_{|r|\lesssim Zq/X}\widehat{K_1}(r\ov{q_0})\frac{1}{q_0^{1/2}}\sumstar_{\alpha(q_0)}\widehat{K_0}((r-\alpha q_1)\ov{q_1})e\left(\frac{\mp\bar{\alpha} n}{q_0}\right) \widehat{V}\left(\frac{rX}{q}\right).
\end{split}\end{equation}
  Here we have made use of the twisted multiplicativity 
\eqref{twisted-multiplicity}.
  Concerning the $\alpha$-sum we have
  \begin{equation}\label{character-in-error}
\frac{1}{q_0^{1/2}}\sumstar_{\alpha({q_0})}\widehat{K_0}((r-\alpha {q_1})\overline{{q_1}})e\left(\frac{\mp\bar{\alpha} n}{{q_0}}\right)
=\frac{1}{q_0^{1/2}}\sum_{x({q_0})}K_0(x)e\left(\frac{r\overline{{q_1}}x}{{q_0}}\right)\mathrm{Kl}_2(\pm nx;{q_0})\ll_{C_0}  1.
\end{equation}
That follows from \cite[Lem. 8.1]{LMS} and our assumption that the sheaf $\mcF$ associated to $K_0$ satisfies the (MO) condition. 

Therefore we obtain 
  \begin{equation}\label{error-c=q0}
\begin{split}
|S(X)_{c=q_0}|\ll \|\widehat{K_1}\|_{\infty}\frac{X^{2+o(1)}}{P{q_0}^{3/2}q^{1/2}}\frac{{q_0}^2}{X}\frac{Zq}{X}\ll \|\widehat{K_1}\|_{\infty} X^{o(1)}\frac{Zq_0^{1/2}q^{1/2}}{P}.
\end{split}\end{equation}

  \subsection{Bounding $S(X)$: Conclusion}  
Plugging the bounds in \eqref{S-off} and \eqref{S-dia} for $S_{\text{off}}(X)$ and $S_{\text{dia}}(X)$ into \eqref{final-object}, we eventually obtain that 
 \begin{equation*}
\begin{split}
S(X)\ll&  \|\widehat{K_1}\|_{\infty}  Z^{1/2}X^{1/2+o(1)} q_0^{1/2}+\|\widehat{K_1}\|_{\infty} ZX^{1/2+o(1)}q^{1/2}{q_0}^{-1/4}\\
&\quad\quad\quad\quad \quad\quad +\|\widehat{K_1}\|_{\infty} ZX^{o(1)}q^{1/2}{q_0}^{1/4}+|S(X)_{c=p}|+|S(X)_{c=q_0}|\\
\ll&\|\widehat{K_1}\|_{\infty}X^{o(1)}\Bigl(  Z^{1/2}X^{1/2} q_0^{1/2}+ ZX^{1/2}q^{1/2}{q_0}^{-1/4}\\
&\quad\quad\quad\quad\quad\quad \quad\quad + Zq^{1/2}{q_0}^{1/4}+\frac{ ZX^{1/2}q^{1/2}}{q_0}+ \frac{Zq_0^{1/2}q^{1/2}}{P}\Bigr)
\end{split}
 \end{equation*}
 upon inserting the bounds \eqref{error-c=p} and \eqref{error-c=q0} for $|S(X)_{c=p}|+|S(X)_{c=q_0}|$.
To satisfy \eqref{initial-assumption}, one needs
\begin{equation*}
p{q_0}\gg X
\end{equation*}
which can be met by taking $P$ ($<q^A$) large enough.

\begin{remark}
(1) The parameter $P$ does not play an essential role in the final bounds (cf. see also \cite{AHLQ}).

\par

(2) Similar results can also be proved using an amplified second moment approach as in \cite{Bykovskii} and \cite{FKM1}.
\end{remark}

\section{Proof of Corollary \ref{mainarith}}

In this section we show how to derive Corollary  \ref{mainarith} from  Theorem \ref{mainthm2}.

 Let $q=q_0q_1$ be a product of two primes as above. We apply the duality principle \cite[Cor. 9.2]{LMS} to the sum
$$S_V(K;X)=\sum_{n|r|\geq 1}\lambda_\varphi(n,r)\lambda_f(n)K(nr^2)V\left(\frac{nr^2}{X}\right)$$
where
$$K(n)=q^{1/2}\delta_{n\equiv a\bmod  q}.$$
Setting $\widecheck X=q^6/X$, this gives (up to negligible error terms) 
$$
q^{-1/2}S_V(K;X)=\frac{X}{q^{7/2}}\sum_{r|n|\geq 1}\lambda_\varphi(r,n)\lambda_f(n)\mathrm{Kl}_6(anr^2;q)\widecheck{V}\left(\frac{r^2n}{\widecheck{X}}\right).$$
By the Chinese Remainder Theorem we have
$$\Kl_6(n;q_0q_1)=\Kl_6(\ov q_1^6n;q_0)\Kl_6(\ov q_0^6n;q_1)$$
and the hyper-Kloosterman sum $\Kl_6(\ov q_1^6n;q_0)$ is the trace function attached to the hyper-Kloosterman sheaf $[\times \ov q_1^6]^*\KL_6$ which is ``good" in the sense of \cite{LMS}. Applying Theorem \ref{mainthm2} we obtain
\begin{equation*}
\begin{split}
\sum_\stacksum{n,|r|\geq 1}{nr^2\equiv a\bmod q}\lambda_\varphi(n,r)\lambda_f(n)V\left(\frac{nr^2}{X}\right)
\ll& X^{o(1)} \frac{X}{q^{7/2}}\left(\widecheck{X}^{3/4}q^{3/5}+\widecheck{X}q^{-1/5}\right)\\
=&X^{o(1)}\left(X^{1/4}q^{8/5}+q^{23/10}\right)\leq (X/q)^{1-\delta}
\end{split}
\end{equation*}
for some $\delta>0$ as long as 
$$q\leq  X^{15/52-\eta}=X^{\theta_6+1/364-\eta}$$ for some $\eta>0$ (here $\theta_6=2/7$ cf. \eqref{thetaddef}).

\section{Proof of Theorem \ref{mainthm2}} 
In this section, we give a detailed proof for Theorem \ref{mainthm2}, following the same approach as in \cite{LMS}. The proof will be very similar to the one presented in \cite{LMS}, the chief difference being in the choice of the parameter when the delta method is applied. Since that paper contains almost all the necessary technical details (corresponding to the case $q_0=q$) we will be brief here (with some simplifications) and refer the readers to \cite{LMS} for the relevant details (see also \cite{Sun-Yu} where the case $K=\chi$ a multiplicative character is treated).

\subsection{First transformations} 
From now on we assume that $Z$ satisfies
$1\leq Z\leq q$ and  $V\in \mcC^\infty_c(\Rr)$  satisfying the bound \eqref{bound-of-V}.
As in \cite{LMS}, we write
\begin{equation}\label{Sstart}
S^t_V(K,X):=\sum_{r,n}\lambda(r,n)\lf(n)K(nr^2)V\left(\frac{nr^2}{X}\right)=\sum_{|r|\leq R}S_{V,r}(K,X/r^2)+X^{o(1)}R^{\theta_3}X/R,
\end{equation}
where
\begin{equation}\label{Sdef}
S_{V,r}(K,X):=\sum_{n=1}^{\infty}\lambda(r,n)\lambda_f(n)K(nr^2)V\left(\frac{n}{X}\right)
\end{equation}
 and $R$ satisfying $R<q$ is some parameter to be determined later (see \eqref{choice-of-R}).
In the above $\theta_3=5/14$ is the Kim--Sarnak bound towards the Ramanujan--Petersson conjecture on $\GL_3$. For each fixed $r$ with $r\leq R$ we can write $S_{V,r}(K,X)$ as 
\begin{equation}
     S_{V,r}(K,X)=\sum_{n=1}^{\infty}\lambda(r,n)\sum_{m=1}^{\infty}\lf(m)K(mr^2)\delta_{n=m}U\left(\frac{n}{X}\right)V\left(\frac{m}X\right).    \label{SVXfirst}
\end{equation}
Here $U$ is a smooth function supported in $(1/100,100)$ and satisfying $U^{(j)}(x)\ll_j 1$ for $j\geq 0$ and $U(x)=1$ for $x\in [1,2]$.

Let $C=(X/q_0)^{1/2}$. We apply a version of the Duke--Friedlander--Iwaniec delta method (\cite{DFI1.5})
\begin{eqnarray*}
\delta_{n=0}&=&\frac{1}{C}\sum_\stacksum{c\leq  2C}{(c,q)=1}\frac{1}{c{q_0}}\sumstar_{u(c{q_0})}e\left(n\frac{u}{c{q_0}}\right)h\left(\frac{c}{C},\frac{n}{C^2{q_0}}\right)\nonumber\\
&+&\frac{1}{C}\sum_\stacksum{c\leq  2C}{(c,q)=1}\frac{1}{c{q_0}}\sumstar_{a(c)}e\left(n\frac{a}{c}\right)h\left(\frac{c}{C},\frac{n}{C^2{q_0}}\right)\label{eqdelta}\\
&+&
\frac{1}{C}\sum_{c\leq  2C/q}\frac{1}{cq{q_0}}\sumstar_{a(cq{q_0})}e\left(n\frac{a}{cq{q_0}}\right)h\left(\frac{cq}{C},\frac{n}{C^2{q_0}}\right)+O_A\left(C^{-A}\right)\nonumber
\end{eqnarray*}
as presented in \cite[Thm. 1]{HB} and \cite[(3.7)]{LMS} (see also \cite[Lem. 2]{Sun-Yu} for a similar version), to the difference $n-m$ in \eqref{SVXfirst} and obtain
\begin{equation}\label{S'sumErr}
S_{V,r}(K,X)=\mathrm{Main}+\mathrm{Err}_1+\mathrm{Err}_2+O_A\left(X^{-A}\right)
\end{equation}
where
\begin{multline}\label{S'sum}
       \mathrm{Main}=\frac{1}{Cq_0}\sum_\stacksum{c\leq  2C}{(c,q)=1}\frac{1}{c}\sumstar_{u(cq_0)}\sum_{n=1}^{\infty}\lambda(r,n)e\left(\frac{un}{cq_0}\right)U\left(\frac{n}{X}\right)\\
        \times\sum_{m=1}^{\infty}\lf(m)K(mr^2)e\left(\frac{-um}{cq_0}\right)V\left(\frac{m}X\right)h\left(\frac{c}{C},\frac{n-m}{C^2q_0}\right),
\end{multline}
and
\begin{multline}\label{eqErr2}
        \mathrm{Err}_1=\frac{1}{Cq_0}\sum_\stacksum{c\leq  2C}{(c,q)=1}\frac{1}{c}\sumstar_{a(c)}\sum_{n=1}^{\infty}\lambda(r,n)e\left(\frac{an}{c}\right)U\left(\frac{n}{X}\right)\\
        \times\sum_{m=1}^{\infty}\lf(m)K(mr^2)e\left(\frac{-am}{c}\right)V\left(\frac{m}X\right)h\left(\frac{c}{C},\frac{n-m}{C^2q_0}\right),
\end{multline}
\begin{multline}\label{eqErr3}
        \mathrm{Err}_2=\frac{1}{Cq{q_0}}\sum_\stacksum{c\leq  2C/q}{(c,q)=1}\frac{1}{c}\sumstar_{a(cqq_0)}\sum_{n=1}^{\infty}\lambda(r,n)e\left(\frac{an}{cqq_0}\right)U\left(\frac{n}{X}\right)\\
        \times\sum_{m=1}^{\infty}\lf(m)K(mr^2)e\left(\frac{-am}{cqq_0}\right)V\left(\frac{m}X\right)h\left(\frac{cq}{C},\frac{n-m}{C^2{q_0}}\right).
\end{multline}

In the following, we focus our analysis on the term $\mathrm{Main}$ in \eqref{S'sum} which is the hardest and is responsible for the final bound. The treatment for the terms $\mathrm{Err}_1,\mathrm{Err}_2$ are almost identical to the one presented in \cite[\S 7.1]{LMS}, and just as in \cite{LMS} their contribution turns out to be smaller as compared to that of $\mathrm{Main}$. As such we completely skip their treatments here and the reader is referred to \cite[\S 7.1]{LMS} for the details.

\subsubsection{Bounding $\mathrm{Main}$} 

To prepare for the application of Voronoi summation formula to the $m$-sum in \eqref{S'sum} we write
$$K(mr^2)=\frac{1}{q^{1/2}}\sum_{b\bmod  q}\widehat{K}(b)e(\frac{-bmr^2}{q})=
\frac{1}{q^{1/2}}\sum_{b\bmod  q}\widehat{K}(b)e(\frac{-bcr^2m}{cq}),$$
where $\widehat{K}(b)=\widehat{K}(b;q)$
denotes the Fourier transform of $K$; see \eqref{FourierTransform-K}.

We find that the term in \eqref{S'sum} can be rewritten as
\begin{multline}\label{S'sum-2}
       \mathrm{Main}=\frac{1}{Cq_0}\sum_\stacksum{c\leq  2C}{(c,q)=1}\frac{1}{c}\sumstar_{u(cq_0)}\sum_{n=1}^{\infty}\lambda(r,n)e\left(\frac{un}{cq_0}\right)U\left(\frac{n}{X}\right)\\
        \times\frac{1}{q^{1/2}}\sum_{b\bmod  q}\what K(b)\sum_{m=1}^{\infty}\lf(m)e\left(\frac{-(bcr^2+u q_1)m}{cq}\right)V\left(\frac{m}X\right)h\left(\frac{c}{C},\frac{n-m}{C^2q_0}\right).
\end{multline}
We can further assume that $(bcr^2+u q_1,cq)=1$, as otherwise we would have $(bcr^2+u q_1,q)\geq q_1$ and the contribution from the latter case can be seen to be much smaller (see \cite[\S 3.1]{LMS}).
        
Under $(bcr^2+u q_1,cq)=1$, we now apply Voronoi summation to the $m$-sum in \eqref{S'sum-2}, to get
\begin{multline}\label{eqafter1stvoronoi}
              \mathrm{Main}=\frac{X}{Cq_0q}\sum_{\pm}\sum_\stacksum{c\leq  2C}{(c,q)=1}\frac{1}{c^2}\frac{1}{q^{1/2}}\sumsumstar_\stacksum{b(q),u(cq_0)}{(bcr^2+u q_1,cq)=1}\what K(b)\\
        \times\sum_{m=1}^{\infty}\ov{\lf(m)}e(\frac{\pm\ov{bcr^2+u q_1}m}{cq})\sum_{n=1}^{\infty}\lambda(r,n)e\left(\frac{un}{cq_0}\right)U\left(\frac{n}{X}\right)\widehat{\mathcal{V}}^{\pm}\left(n,\frac{mX}{c^2q^2}\right)+O\left(X^{-A}\right),
        \end{multline}
where $\widehat{\mathcal{V}}^{\pm}\left(n,y\right)$ is given as in \cite[(3.16)]{LMS}.

Next, we further apply Voronoi summation (\cite[Prop. 2.2]{LMS}) to the $n$-sum above to obtain a sum of the form
\begin{equation}
\label{eqaftervoronoi}
\frac{Xq_0^{1/2}}{Cq}\sum_{\pm\pm}\sum_\stacksum{c\leq  2C}{(c,q)=1}\frac{1}{c}
\sum_{n_1|rcq_0}\sum_{m,n}\ov{\lf(m)}\frac{\lambda(n,n_1)}{nn_1}\mathcal{C}(m,n;\frac{rcq_0}{n_1})\mcW_{\pm\pm}(\frac{m}{c^2q^2/X},\frac{n_1^2n}{c^3q^3_0r/X})    
\end{equation}	

where
$$\mathcal{C}(m,n;\frac{rcq_0}{n_1})=\frac{1}{(qq_0)^{1/2}}\sumsumstar_\stacksum{b(q),u(cq_0)}{(bcr^2+u q_1,cq)=1}\what K(b)e(\frac{\pm\ov{bcr^2+u q_1}m}{cq})S(r\ov u,\pm n;\frac{rcq_0}{n_1})$$
and
\begin{equation*}\label{double-weight-function}
\begin{split}
\mcW_{\pm\pm}(y,z)&=\int_{\mathbb{R}}V(x)\mathcal{W}_{x,\pm}(z)\mcJ^\pm_{f}(4\pi \sqrt{xy}){\mathrm{d}x}\\
=&z\int_{\mathbb{R}}V(x)\bigl(\int_{0}^{\infty}W_x(\xi)J_{\vphi,\pm}(z\xi){\mathrm{d}\xi}\bigr)\, \mcJ^\pm_{f}(4\pi \sqrt{xy}){\mathrm{d}x};
  \end{split}\end{equation*}
see \cite[(3.19)]{LMS}. Here $$W_x(\xi):=U(\xi)h\left(\frac{c}{C},\frac{X(\xi-x)}{C^2{q_0}}\right).$$
In particular, for $z\ll 1$ we have 
\begin{equation}\label{Jacquet-Shalika}
\mcW_{\pm\pm}(y,z)\ll z^{1/2}.
\end{equation}

\subsection{The case $q_0\not| n_1$}
For the sum in \eqref{eqaftervoronoi}, we further split it into two subsums according to $(n_1,q_0)=1$ or not, and write 
$$ \mathrm{Main}= \mathrm{Main}_{00}+ \mathrm{Err}_4+O\left(X^{-A}\right),$$
where
\begin{multline}
\mathrm{Main}_{00}=\frac{Xq_0^{1/2}}{Cq}\sum_{\pm\pm}\sum_\stacksum{c\leq  2C}{(c,q)=1}\frac{1}{c}\sum_\stacksum{n_1|rc}{(n_1,q_0)=1}\sum_{m,n}\ov{\lf(m)}\\
\times\frac{\lambda(n,n_1)}{nn_1}\mathcal{C}(m,n;\frac{rcq_0}{n_1}) \mcW_{\pm\pm}(\frac{m}{c^2q^2/X},\frac{n_1^2n}{c^3q^3_0r/X})\label{n1qcoprime}
\end{multline}
and $\mathrm{Err}_4$ corresponds to the complementary sum where $q_0|n_1$, whose contribution is given in \eqref{qdivn1final} (see \cite[\S  6.3]{LMS}).

We have (since $(c,q)=1$)
$$e(\frac{\pm\ov{bcr^2+u q_1}m}{cq})=e(\frac{\pm\ov{bcr^2+u q_1}\ov{c}m}{q})e(\frac{\pm\ov{uq_1}\ov qm}{c})$$
$$S(r\ov u,\pm n;\frac{rcq_0}{n_1})=S(\ov{c}n_1\ov u,\pm \ov{rc}n_1n;{q_0})
S(\ov{q_0}r\ov u,\pm \ov{q_0}n;{rc}/n_1).$$
Therefore, the $(b,u)$-sum in $\mathcal{C}(m,n;\frac{rcq_0}{n_1})$ splits into a product of two sums of respective moduli $rc/n_1$ and $q_0$. Accordingly, we write
$$\mathcal{C}(m,n;\frac{rcq_0}{n_1})=M_{n_1,r}(m,n;rc)N_{\ov{cr}}(m,n)$$
where
\begin{equation*}\label{Mrcsum}
M_{n_1,r}(m,n;rc):=\sumstar_{u(c)}e(\frac{\pm\ov{uq^2_1}\ov{q_0}m}{c})S(\ov{q_0}r\ov u,\pm \ov{q_0}n;{rc}/n_1),
\end{equation*}
and 
\begin{multline}
N_{\ov{cr}}(m,n):=\frac{1}{(qq_0)^{1/2}}\sumsumstar_\stacksum{b(q),u(q_0)}{(bcr^2+u q_1,q)=1}\what K(b)e(\frac{\pm\ov{bcr^2+u q_1}\ov{c}m}{q})S(\ov{c}n_1\ov u,\pm \ov{rc}n_1n;{q_0})\\=\frac{1}{q^{1/2}}\sumsumstar_\stacksum{b(q),u(q_0)}{(b+uq_1,q)=1}\what K(b)e(\frac{\pm\ov c^2\ov r^2 m\ov{b+uq_1}}{q})\Kl_2(\pm \ov{c}^3\ov{r}^3 n_1^2n\ov u;{q_0})\\
=	\sumstar_{u(q_0)}L_{\pm\ov c^2\ov r^2m,1}(uq_1;q)\Kl_2(\pm \ov{c}^3\ov{r}^3 n_1^2n\ov u;{q_0})
\label{Nq0sum}
\end{multline}
with
\begin{equation}\label{definition-L}
L_{\alpha,\beta}(u;q):=\frac{1}{q^{1/2}}\sum_\stacksum{b(q)}{(b+\beta u,q)=1}\what K(b)e\left(\frac{\alpha\, \ov{b+\beta u}}{q}\right).
\end{equation}

From these notations we find that the sum $\mathrm{Main}_{00}$ in \eqref{n1qcoprime} is equal to
\begin{multline*}\label{after-voronoi}
\mathrm{Main}_{00}=\frac{Xq_0^{1/2}}{Cq}\sum_{\pm\pm}\sum_\stacksum{c\leq  2C}{(c,q)=1}\frac{1}{c}\sum_\stacksum{n_1|rc}{(n_1,q_0)=1}\sum_{m,n}\ov{\lf(m)}\frac{\lambda(n,n_1)}{nn_1}\times\\
M_{n_1,r}(m,n;rc)N_{\ov{cr}}(m,n) \mcW_{\pm\pm}(\frac{m}{c^2q^2/X},\frac{n_1^2n}{c^3q^3_0r/X}).
\end{multline*}

We break the $c$-sum into $O(\log X)$ many dyadic intervals with $c\sim C'$, where $C'$ satisfies
$$ C'\leq 2C=2(X/{q_0})^{1/2}.$$ By \cite[Lem. 3.3]{LMS}  we know that the $(m,n)$-sum can be truncated at
\begin{equation}\label{truncation-of-n}
m\lesssim M=Z^2\frac{{C'}^2q^{2}}{X}, nn^2_1\lesssim \frac{{C'}^3q^{3}_0r}{X}.
\end{equation}

For each fixed $n_1$ we break the $n$-sum into $O(\log q)$ dyadic intervals $n\sim N/n^2_1$ with $N$ satisfying $$N\lesssim \frac{{C'}^3q^{3}_0r}{X}.$$

Now for each $c\sim C'$ and $nn^2_1\sim N$ we will evaluate the truncated version of $\mathrm{Main}_{00}$.

\subsubsection{Cauchy--Schwarz}\label{CSsec}
We now factor $c=c_1c_2$ with $$c_1\leq C',\ {n_1|rc_1},\ c_1|(n_1r)^\infty\hbox{ and }(c_2,n_1r)=1.$$
Then we apply Cauchy--Schwarz inequality and the Rankin--Selberg estimate to bound  the sum $\mathrm{Main}_{00}$ as follows (for the various choices of $\pm,\pm$)
\begin{equation}\label{main-after-cauchy}
\mathrm{Main}_{00}\ll  X^{o(1)}\frac{Xq_0^{1/2}}{CqC'}\sup_{N\lesssim \frac{{C'}^3q^{3}_0r}{X}}\frac{1}{N^{1/2}}B(N)^{1/2}
\end{equation}
with
\begin{multline*}
B(N):=\sumsum_\stacksum{c_1,nn_1^2\approx N}{(n_1,q_0)=1}n_1\biggl| \sum_{m\leq M}\ov{\lf(m)}\\
\times
\sum_\stacksum{c_2\sim C'/c_1}{(c_2,q)=1}M_{n_1,r}(m,n;rc_1c_2)N_{\ov {c_1c_2r}}(m,n)\mcW_{\pm\pm}(\frac{m}{c_1^2c_2^2q^2/X},\frac{n_1^2n}{c_1^3c_2^3q^3_0r/X}) \biggr|^2 U\left(\frac{n}{N/n^2_1}\right).
\end{multline*}
Here $U$ is a smooth function with compact support contained in $(0,\infty)$ satisfying $U^{(j)}(x)\ll_j 1$ for $j\geq 0$.

After opening the square, the factor $B(N)$ equals
\begin{gather}\label{eqprepoisson}
B(N)=\sumsum_{c_1,n_1}n_1\sumsum_{m,m'}\ov{\lf(m)}\lf(m')\sumsum_{c_2,c_2'}\times\\
\sum_{n\geq 1} M_{n_1,r}(m,n;rc_1c_2)\ov{M_{n_1,r}(m',n;rc_1c'_2)}\\
\times N_{\ov {c_1c_2r}}(m,n)\ov{N_{\ov {c_1c'_2r}}(m',n)}\mcW\left(\frac{n}{N/n^2_1}\right),	\nonumber
\end{gather}
 where 
 \begin{equation*}\label{weight-before-poisson}
 \mcW\left(\frac{n}{N/n^2_1}\right)=U\left(\frac{n}{N/n^2_1}\right)\mcW_{\pm\pm}(\frac{m}{c^2_1c^2_2q^2/X},\frac{n_1^2n}{c^3_1c^3_2q^3_0r/X})\ov{\mcW_{\pm\pm}}(\frac{m'}{c^2_1{c'_2}^2q^2/X},\frac{n_1^2n}{c^3_1{c'_2}^3q^3_0r/X}).
 \end{equation*}

We apply Poisson formula to the $n$-variable keeping in mind that 
\begin{equation*}\label{eq-MMNN} 
n\mapsto M_{n_1,r}(m,n;rc_1c_2)\ov{M_{n_1,r}(m',n;rc_1c'_2)}N_{\ov {c_1c_2r}}(m,n)\ov{N_{\ov {c_1c'_2r}}(m',n)}
\end{equation*}
 is periodic of period $q_0rc_1c_2c_2'/n_1:=q_0k$, and see that \eqref{eqprepoisson} equals
\begin{equation}\label{eqpostpoisson}
\sumsum_{c_1,n_1}n_1\sumsum_{m,m'}\ov{\lf(m)}\lf(m')\sumsum_{c_2,c_2'}\frac{N}{n^2_1\sqrt{q_0k}}\sum_{n\in\Zz}\mathrm{FT}(n,m,m';q_0k)\what \mcW(n/N^*),
\end{equation}
 where (after inserting \eqref{Nq0sum})
       \begin{multline*}
\mathrm{FT}(n,m,m';q_0k)=\sumsum_{u,u'\bmod  {q_0}}
L_{\pm\ov c^2\ov r^2m,1}(uq_1;q)\ov{L_{\pm\ov {c'}^2\ov r^2m',1}(u'q_1;q)}\times\\
\frac{1}{\sqrt{q_0k}}\sum_{v\bmod  {q_0k}}\Kl_2(\pm \ov{c}^3\ov{r}^3 n_1^2v\ov u;{q_0})
\ov{\Kl_2(\pm \ov{c'}^3\ov{r}^3 n_1^2v\ov {u'};{q_0})} M_{n_1,r}(m,v;rc)\ov{M_{n_1,r}(m',v;rc')}\, e\left(\frac{nv }{q_0k}\right)
\end{multline*} 
(with $c=c_1c_2$ and $c'=c_1c'_2$) and
     \begin{equation}\label{defN*}
       N^*:=q_0kn_1^2/ N.
       \end{equation}
    Here due to the rapid decay of the weight function $\what \mcW(y)$ when $|y|\gg q^{\varepsilon}$, we can truncate the dual $n$-sum in \eqref{eqpostpoisson} at $|n|\ll q^{\varepsilon} N^*$. By using the estimate \eqref{Jacquet-Shalika} and the truncation in \eqref{truncation-of-n}, for $|n|\ll q^{\varepsilon} N^*$, we readily have the bound (see for instance \cite[(27)]{Lin})
    \begin{equation}\label{scalar-of-W}
    \what \mcW(n/N^*)\ll  \frac{N}{{C'}^3q^{3}_0r/X}.
    \end{equation}

\subsubsection{Computation of $\mathrm{FT}$}\label{computation-fourier-transform}
Recall $k=rc_1c_2c'_2/n_1$. We have $(q_0,k)=1$ and we split the sum $\mathrm{FT}(n,m,m';q_0k)$ in the above as a product of sums $\FT(n)$ and $\FT(n;k)$ of respective moduli $q_0$ and $k$ (to simplify notations we do not display the dependency in $m,m'$ in these expressions).
  \subsubsection{The $k$-sum} The $k$-sum equals
$$\FT(n;k):=\frac{1}{\sqrt{k}}\sum_{v(k)}M_{n_1,r}(m,v;rc_1c_2)\ov{M_{n_1,r}(m',v;rc_1c'_2)}\, e\left(\frac{ nv \ov{q_0}}k\right).$$
The following bounds can be found in \cite[Lem. 4.4]{LMS}.
\begin{lemma}\label{bound-k-part}
 We have the following estimates
\begin{equation*}
\FT(0;k)\ll \sqrt{k}rc_1c_2\mathop{\sum_{d|c_1c_2}\sum_{d'|c_1c_2}}_{(d,d')|(m  -m')} (d ,d'),
\end{equation*}
and
\begin{equation*}
\begin{split}
\FT(n;k)\ll \sqrt{k}\sum_{d_1|c_1}d_1\sum_{d'_1|c_1}d'_1
\mathop{\sumstar_\stacksum{x_1(rc_1/n_1)}
{q^2_1  n_1x_1\equiv \mp m\bmod {d_1}}}\sumsum_\stacksum{d_2|(c_2,q^2_1 n_1c'_2+nm)}{d'_2|(c'_2,q^2_1 n_1c_2+nm')}d_2d'_2.
\end{split}\end{equation*}
\end{lemma}

\subsubsection{The $q_0$-sum}    \label{secq-sum}   
       The $q_0$-sum equals
\begin{multline}\label{eq-qsum2}
\FT(n)=\frac{1}{\sqrt {q_0}}\sumsum_{u,u'\bmod  {q_0}}
L_{\pm\ov c^2\ov r^2m,1}(uq_1;q)\ov{L_{\pm\ov {c'}^2\ov r^2m',1}(u'q_1;q)}\\
\times\sum_{v\bmod  {q_0}}\Kl_2(\pm \ov{c}^3\ov{r}^3 n_1^2v\ov u;{q_0})\ov{\Kl_2(\pm \ov{c'}^3\ov{r}^3 n_1^2v\ov {u'};{q_0})}\, e\left(\frac{\ov kv n}{q_0}\right).
\end{multline}

As we recall, by making use of the factorisation $q=q_0q_1$, we have the twisted multiplicativity \eqref{twisted-multiplicity}
$$    \widehat{K}(b;q)=\widehat{K_0}(\overline{{q_1}}b) \widehat{K_1}(\overline{{q_0}}b)$$
and then (recalling \eqref{definition-L})
$$L_{\alpha,\beta}(uq_1;q)=L_{\alpha\ov{q_1}^2,\beta}(u)L_{\alpha\ov{q_0}^2,\beta}(0),$$
where
$$L_{\alpha,\beta}(u):=\frac{1}{\sqrt{q_0}}\sum_{b(q_0)}\widehat{K_0}(b)e\left(\frac{\alpha\overline{(b+\beta u)}}{q_0}\right), L_{\alpha,\beta}(u):=\frac{1}{\sqrt{q_1}}\sum_{b(q_1)}\widehat{K_1}(b)e\left(\frac{\alpha\overline{(b+\beta u)}}{q_1}\right).$$

With this factorisation we can write \eqref{eq-qsum2} as
\begin{multline*}
\FT(n)=L_{\pm\ov c^2\ov r^2\ov{q_0}^2m,1}(0)\ov{L_{\pm\ov {c'}^2\ov r^2\ov{q_0}^2m',1}(0)}\times\\
\frac{1}{\sqrt {q_0}}\sumsum_{u,u'\bmod  {q_0}}
L_{\pm\ov c^2\ov r^2\ov{q_1}^2m,1}(u)\ov{L_{\pm\ov {c'}^2\ov r^2\ov{q_1}^2m',1 }(u')}\times\\
\sum_{v\bmod  {q_0}}\Kl_2(\pm \ov{c}^3\ov{r}^3 n_1^2v\ov u;{q_0})\Kl_2(\pm \ov{c'}^3\ov{r}^3 n_1^2v\ov {u'};{q_0})\, e\left(\frac{\ov kv n}{q_0}\right).
\end{multline*}

According to the calculations in \cite[\S 4.2.2]{LMS}, we can further express
$$
\FT(n)=L_{\pm\ov c^2\ov r^2\ov{q_0}^2m,1}(0)\ov{L_{\pm\ov {c'}^2\ov r^2\ov{q_0}^2m',1}(0)}\times
\sqrt{q_0}\sum_{v\bmod  {q_0}}Z(v)\ov{Z'(v-\delta)},$$
where
\begin{equation*}\label{Zcompute}
Z(v)=Z_{\alpha,\beta,\gamma}(v):=\frac{1}{\sqrt{q_0}}\sum_{x\in\Fqqt}\Kl_2(\beta\gamma x)K_0(xv)\Kl_2(\alpha xv),	
\end{equation*}
and $Z'(v)$ is defined likewise with the parameter $(\alpha,\beta,\gamma)$ being replaced by $(\alpha',\beta',\gamma')$.

The following choices of values of the parameters correspond to our initial problem:
\begin{gather}\alpha=\pm \ov c^2\ov r^2\ov{q_1}^2m,\  \alpha'=\pm \ov {c'}^2\ov r^2\ov{q_1}^2m',\ \beta=\beta'=1 \nonumber\\
\gamma=\pm \ov{c}^3\ov{r}^3 n_1^2,\ \gamma'=\pm \ov{c'}^3\ov{r}^3 n_1^2,\ \delta=\ov k n,\label{eq-actual}
\end{gather}

We recall the following result proven in \cite[\S 8]{LMS}.
\begin{proposition}\label{sqrootcancel} 
Let $T_\mcF(\Fq)$ be the subgroup of $\Fqt$ defined by
$$ T_\mcF(\Fq)=\{\lambda\in\Fqt,\ [\times\lambda]^*\mcF\hbox{ is geometrically isomorphic to }\mcF\}.
$$
 Assuming that the sheaf $\mcF$ is good, then for any $\alpha,\beta,\alpha',\beta',\gamma,\gamma',\delta\in\Fqt$, we have
$$\sum_{v}Z(v)\ov{Z'(v-\delta)}=O(q^{1/2}).$$
If $\delta=0$ the above bound holds unless
$$\alpha/\alpha'=\beta\gamma/\beta'\gamma'\in T_\mcF(\Fq)$$
in which case
$$\sum_{v}Z(v)\ov{Z'(v)}=c_\mcF(\alpha/\alpha')q+O(q^{1/2})$$
for $c_\mcF(\alpha/\alpha')$ some complex number of modulus $1$.
Here the implicit constants depend only on $C(\mcF)$.
\end{proposition}

Returning to our original sum, applying \cite[Lem. 8.1]{LMS} to 
$$L_{\alpha,1}(0)=\sum_{x\bmod  {q_1}}K_1(x)\Kl_2(\alpha x)$$
and Proposition \ref{sqrootcancel} (with $q=q_0$), we see that the $\FT(n)$ in \eqref{eq-qsum2} is  $O(q_0)$ unless $\delta=0$ and
$${c'}^2m/c^2m'={c'}^3/c^3\in T_\mcF({\mathbf{F}_{q_0}})$$ in which case 
\eqref{eq-qsum2} equals $C({c'}^2m/c^2m'){q_0}^{3/2}+O(q_0)$ with $|C({c'}^2m/c^2m')|=1.$

\subsection{Contribution of the $n=0$ frequency}\label{contribution-of-zero}

In this section we bound the contribution  to \eqref{eqpostpoisson} from the frequency $n=0$ in \eqref{eqprepoisson}. By \cite[(4.15)]{LMS}, we then have
\begin{equation}\label{cequal}
c_2=c_2',\ c=c',\ k=rc_1c_2^2/n_1.	
\end{equation}

We use the case $\delta=0$ of Proposition \ref{sqrootcancel}: by \eqref{cequal} and \eqref{eq-actual} we have that \eqref{eq-qsum2} is  $O(q_0)$ unless we have the congruence modulo $q_0$
$$m/m'=1\in T_\mcF({\mathbf{F}_{q_0}})$$ in which case 
\eqref{eq-qsum2} equals $C(m/m'){q_0}^{3/2}+O(q_0)$ with $|C(m/m')|=1.$
 
According to the calculations in \cite[\S 5]{LMS}, the contribution of the $n=0$ frequency to \eqref{eqpostpoisson} is bounded by
$$X^{o(1)}
\frac{rN{C'}^2M}{q_0^{1/2}}(C'{q_0}^{3/2}+q_0{M}).$$

Taking the square root of this term and multiplying it by $\frac{Xq_0^{1/2}}{CqC'}\frac{1}{N^{1/2}}$ where $N\ll \frac{{C'}^3q^{3}_0r}{X}$, we see that the contribution of these terms to \eqref{main-after-cauchy} and therefore to \eqref{n1qcoprime} is bounded by
\begin{multline}\label{n=0bound}
X^{o(1)}\frac{r^{1/2}XM^{1/2}{q_0}^{1/4}}{Cq}\left({C'}^{1/2}{q_0}^{3/4}+q_0^{1/2}M^{1/2}\right)\\
\ll X^{o(1)}r^{1/2}\left(ZX^{3/4}{q_0}^{3/4}+Z^2 X^{1/2}q{q_0}^{1/4}\right).
\end{multline}

\subsection{Contribution from the $n\not=0$ frequencies} \label{secn=0}

Recall from \eqref{eqpostpoisson} that
\begin{multline*}
B_{n\neq 0}(N)=\frac{N}{\sqrt{{q_0}}}\sumsum_{c_1,n_1}\frac1{n_1}\sum_{m\leq M}\ov{\lf(m)}\sum_{m'\leq M}{\lf(m')}\sum_{c_2\sim C'/c_1}\sum_{c_2'\sim C'/c_1}\\
\times\frac{1}{\sqrt{k}}\sum_{n\not=0}\mathrm{FT}(n;{q_0})\mathrm{FT}(n;k)\what \mcW(n/N^*).
\end{multline*}

We consider two cases: $n\not=0\bmod  {q_0}$ and $n\equiv0\bmod  {q_0}$.

The following two bounds (by also multiplying the bound \eqref{scalar-of-W} for $\what \mcW(n/N^*)$)
\begin{equation*}
\begin{split}
B_{n\neq 0 \bmod  {q_0}}(N)\ll &X^{o(1)}r^2{C'}^5{q_0}^{3/2}M\left(1+\frac{M}{{C'}}\right)\cdot    \frac{N}{{C'}^3q^{3}_0r/X}\\
B_{{q_0}|n,n\not=0}(N)\ll& X^{o(1)}r^2{C'}^5{q_0}^{3/2}M
\left(1+\frac{M}{{C'}}\right)\frac1{q_0^{1/2}}\cdot    \frac{N}{{C'}^3q^{3}_0r/X}
\end{split}\end{equation*}
were proven in \cite[\S 6.1]{LMS} and \cite[\S 6.2]{LMS} respectively.

The non-zero frequencies contribution to \eqref{main-after-cauchy} and hence to \eqref{n1qcoprime} is bounded by 
\begin{gather}\nonumber
X^{o(1)} \frac{Xq_0^{1/2}}{CqC'} \sup_{N\lesssim \frac{{C'}^3q^{3}_0r}{X}}\frac{1}{N^{1/2}} \left(B_{n\neq 0 \bmod  {q_0}}(N)+B_{{q_0}|n,n\not=0}(N)\right)^{1/2}\\
\ll X^{o(1)} \frac{Xq_0^{1/2}}{CqC'} \sup_{N\lesssim \frac{{C'}^3q^{3}_0r}{X}}\frac{1}{N^{1/2}}   \left(r^2{C'}^5{q_0}^{3/2}M\left(1+\frac{M}{{C'}}\right)\cdot    \frac{N}{{C'}^3q^{3}_0r/X}\right)^{1/2}\nonumber\\
 \ll X^{o(1)}\frac{X^{3/2}}{r^{1/2}C{C'}^{5/2}qq_0}  \left(r^2{C'}^5{q_0}^{3/2}M\left(1+\frac{M}{{C'}}\right)\right)^{1/2}\nonumber\\
\ll X^{o(1)} r^{1/2}\left(\frac{ZX}{{q_0}^{1/4}}+ \frac{Z^{2}X^{3/4}q}{q_0^{1/2}}\right),\label{nonzerocomb}
\end{gather}
since by \eqref{truncation-of-n} $M=Z^2\frac{{C'}^2q^{2}}{X}$ and $C'\leq 2C=2(X/q_0)^{1/2}$.

\subsection{Bounding $S^t_V(K,X)$: the final steps}\label{q-divide-n1}
Let us recall that the sum $\mathrm{Main}$ in \eqref{eqaftervoronoi} was split into two subsums depending on whether $(n_1,q_0)=1$ or not. 

By \eqref{n=0bound} and \eqref{nonzerocomb} the first subsum \eqref{n1qcoprime} is bounded by
\begin{equation}\label{n1qcoprimefinal}
\ll X^{o(1)} r^{1/2}\Big(ZX^{3/4}{q_0}^{3/4}+Z^2X^{1/2}q{q_0}^{1/4}
+\frac{ZX}{{q_0}^{1/4}}+ \frac{Z^{2}X^{3/4}q}{q_0^{1/2}}\Big).	
\end{equation}
According to \cite[(6.5)]{LMS},  the complement sum, $\mathrm{Err}_4$ (when $q_0|n_1$), is bounded by \begin{equation}\label{qdivn1final}
\ll X^{o(1)} {q_0}^{\theta_3}rZ^2\frac{X}{q_0}.	
\end{equation}
Combining this bound with \eqref{n1qcoprimefinal} we see that the sum $\mathrm{Main}$ in \eqref{S'sum} and hence the sum $S_{V,r}(K,X)$ in \eqref{Sdef} is bounded as follows
\begin{multline*}
S_{V,r}(K,X)\ll X^{o(1)}rZ^2\frac{X}{{q_0}^{1-\theta_3}}+	X^{o(1)}r^{1/2}\Big(ZX^{3/4}{q_0}^{3/4}
+Z^2X^{1/2}q{q_0}^{1/4}
+\frac{ZX}{{q_0}^{1/4}}+ \frac{Z^{2}X^{3/4}q}{q_0^{1/2}}\Big). 	
\end{multline*}
Replacing $X$ by $X/r^2$ and averaging this bound over $|r|\leq R$ we obtain
\begin{multline*}
\sum_{|r|\leq R}S_{V,r}(K,X/r^2)\ll X^{o(1)}Z^2\frac{X}{{q_0}^{1-\theta_3}}+	X^{o(1)}\Big(ZX^{3/4}{q_0}^{3/4}
+Z^2R^{1/2}X^{1/2}q{q_0}^{1/4}
+\frac{ZX}{{q_0}^{1/4}}+ \frac{Z^{2}X^{3/4}q}{q_0^{1/2}}\Big).
\end{multline*}
Substituting this into \eqref{Sstart}, we get
\begin{multline*}
S^t_V(K,X)\ll X^{o(1)}Z^2\frac{X}{{q_0}^{1-\theta_3}}+	X^{o(1)}\Big(ZX^{3/4}{q_0}^{3/4}
+Z^2R^{1/2}X^{1/2}q{q_0}^{1/4}
+\frac{ZX}{{q_0}^{1/4}}+ \frac{Z^{2}X^{3/4}q}{q_0^{1/2}}+R^{\theta_3-1}X\Big).
\end{multline*}
By choosing
\begin{equation}\label{choice-of-R}
R=\left(\frac{X}{Z^4q^2q_0^{1/2}}\right)^{\frac{1}{3-2\theta_3}}
\end{equation}
to equate the second and the fifth terms inside the parentheses, we see that
\begin{equation*}
\begin{split}
S^t_V(K,X)\ll& X^{o(1)}Z^2\frac{X}{{q_0}^{1-\theta_3}}+	X^{o(1)}\Big(ZX^{3/4}{q_0}^{3/4}
+Z^{\frac{4(1-\theta_3)}{3-2\theta_3}}
X^{\frac{2-\theta_3}{3-2\theta_3}}(q^2q_0^{1/2})^{\frac{1-\theta_3}{3-2\theta_3}}+\frac{ZX}{{q_0}^{1/4}}+ \frac{Z^{2}X^{3/4}q}{q_0^{1/2}}\Big)\\
\ll& X^{o(1)}Z^2\left(X^{3/4}{q_0}^{3/4}
+
X^{\frac{2-\theta_3}{3-2\theta_3}}(q^2q_0^{1/2})^{\frac{1-\theta_3}{3-2\theta_3}}+\frac{X}{{q_0}^{1/4}}+ \frac{X^{3/4}q}{q_0^{1/2}}\right).
\end{split}
\end{equation*}
In the above, to guarantee $R\geq 1$, we need to assume
\begin{equation}\label{constraint-on-R}
X\geq Z^4q^2q_0^{1/2}.
\end{equation}

\subsection*{Acknowledgements}  Part of this work was inspired by \cite{AHLS}. The first named author would like to thank K. Aggarwal, R. Holowinsky, and Q. Sun for generously sharing their ideas while working on that project. We thank P. Sharma for a helpful comment on an earlier version of this manuscript.

As mentioned in the beginning of this paper, this work is dedicated to the memory of Chandra Sekhar Raju who passed away one year ago, in India, surrounded by his family. Before his premature passing, Chandra was a postdoc in the TAN group at EPFL and before that, a PhD student at Stanford under the supervision of K. Soundararajan. Chandra was not only a gentle soul but also a wonderful colleague, full of ideas, with whom we had countless discussions\footnote{In his thesis \cite{subGL2GL2} Chandra showed, using the $\delta$-symbol method, how to pass the convexity range for character twists of $\GL_2\times\GL_2$ Rankin--Selberg $L$-functions; a lot of our discussions had been revolved around our tentative attempt (so far unsuccessful) to replace Dirichlet characters by more general trace functions.}; several of his ideas are to be found in this paper. It is therefore with both sadness and gratitude that we dedicate this paper to the memory of him.

\end{document}